\documentclass[12pt]{amsart}

\pdfoutput=1

\usepackage{url}
\usepackage{amsmath, amsthm, amssymb, color, bm}
\usepackage{graphicx}
\usepackage{float}
\usepackage{caption}
\usepackage{subcaption}
\usepackage{mathtools}
\usepackage{faktor}
\usepackage{enumerate}
\usepackage{indentfirst}
\usepackage{verbatim}
\usepackage{chemfig}
\usepackage{chemnum}
\usepackage{tikz}
\usepackage{tikz-cd}
\usepackage{caption}
\usepackage{enumitem}
\usepackage{makecell}
\usepackage{array}

\usepackage{geometry}
\usepackage[unicode]{hyperref}
\hypersetup{colorlinks=true,  linkcolor=blue,  filecolor=blue,  citecolor=blue, urlcolor=blue}
\usepackage[capitalize,noabbrev]{cleveref}

\DeclareMathOperator{\val}{val}

\DeclareMathOperator{\GL}{GL}	
\DeclareMathOperator{\SL}{SL}

\DeclareMathOperator{\Ker}{Ker}
\DeclareMathOperator{\diag}{diag}
\DeclareMathOperator{\Nr}{\mathbf{Nr}}
\DeclareMathOperator{\N}{\mathbf{N}}

\DeclareMathOperator{\J}{\mathbf{J}}
\DeclareMathOperator{\D}{\mathbf{D}}
\DeclareMathOperator{\Char}{char}
\DeclareMathOperator{\Gr}{Gr}

\DeclareMathOperator{\SO}{\mathrm{SO}}
\DeclareMathOperator{\Sp}{\mathrm{Sp}}

\newcommand*\A{\mathbb{A}}
\newcommand*\B{\mathbb{B}}
\newcommand*\Z{\mathbb{Z}}
\newcommand*\Q{\mathbb{Q}}
\newcommand*\R{\mathbb{R}}

\newcommand*\E{\mathbb{E}}
\newcommand*\F{\mathbb{F}}
\renewcommand*\P{\mathbb{P}}

\newcommand*\Qp{\Q_p}

\newcommand*\Acal{\mathcal{A}}
\newcommand*\Bcal{\mathcal{B}}
\newcommand*\Ecal{\mathcal{E}}
\newcommand*\Hcal{\mathcal{H}}
\newcommand*\Lcal{\mathcal{L}}
\newcommand*\Mcal{\mathcal{M}}

\newcommand*\Ocal{\mathcal{O}}

\newcommand*\Vcal{\mathcal{V}}
\newcommand*\Xcal{\mathcal{X}}

\newcommand{\norm}[1]{\left\lVert#1\right\rVert}

\newcommand*\varhrulefill[1][2pt]{\leavevmode\leaders\hrule height#1\hfill\kern0pt}

\newcommand{\rvline}{\hspace*{-\arraycolsep}\vline\hspace*{-\arraycolsep}}

\newtheorem{theorem}{Theorem}[section]
\newtheorem{proposition}[theorem]{Proposition}
\newtheorem{lemma}[theorem]{Lemma}
\newtheorem{corollary}[theorem]{Corollary}
\theoremstyle{definition}
\newtheorem{definition}[theorem]{Definition}

\newtheorem{remark}[theorem]{Remark}

\newtheorem{example}[theorem]{Example}

\allowdisplaybreaks

\title{Sampling from $p$-adic algebraic manifolds}

\author{Yassine El Maazouz and Enis Kaya}

\date{\today}


\definecolor{lightyellow}{RGB}{255, 255, 197}

\begin{document}

\begin{abstract}
    We present a method for sampling points from an algebraic manifold, either affine or projective, defined over a local field, with a prescribed probability distribution. Inspired by the work of Breiding and Marigliano \cite{RandomPoints} on sampling real algebraic manifolds, our approach leverages slicing the given variety with random linear spaces of complementary dimension. We also provide an implementation of this sampling technique and demonstrate its applicability to various contexts, including sampling from linear $p$-adic algebraic groups, abelian varieties, and modular curves.

\end{abstract}

\maketitle

\section{Introduction}\label{Sec:1}
    
     For a prime number $p$, the field of $p$-adic numbers $\Q_p$ is obtained by completing the field of rational numbers $\Q$ with respect to the $p$-adic absolute value $|\cdot|_p$. These fields, introduced by Kurt Hensel \cite{Hensel05}, quickly became indispensable in number theory \cite{Hasse1,Hasse2,Mikowski,NeukirchBook,SerreLocalField}. In recent decades, $p$-adic methods have found significant applications in areas such as symbolic computation \cite{LLL,Bostan2005,Kedlaya,Caruso2017} and mathematical physics \cite{GHJMPSST,pAdicAdSCFT,HMPS,padicPhys1}. Moreover, there has been a growing interest in probabilistic and statistical questions in the $p$-adic setting, starting with the pioneering work of Evans \cite{evans2001local}, and extending to more recent contributions \cite{manssour2020probabilistic,BhargavaCremonaFisherGajovic21,BKL22,caruso2021zeroes,maazouz2021gaussian,maazouz2019statistics,evansBMLocalField,evansNoise95,ElementaryDivisors,kulkarni2021integGeo}.
    
    In this paper, inspired by the work of Breiding and Marigliano \cite{RandomPoints}, we address the problem of sampling from a $p$-adic algebraic variety according to a prescribed probability distribution. We do so by intersecting the manifold in question with random linear spaces of complementary dimension. Our results extend the work in \cite{RandomPoints} to the $p$-adic setting and align with $p$-adic analogues of the so-called \emph{Crofton formulas} in integral geometry;  see \cite{BKL22,kulkarni2021integGeo}. These formulas express the volume of a manifold as the expected number of its intersection points with random linear spaces of complementary dimension.

    One of the key motivations for this work stems from the fields of random and numerical algebraic geometry, which have attracted considerable attention in recent decades. For example, in enumerative geometry, the number of solutions of a \emph{generic} zero dimensional system of polynomial equations is usually ill-defined when working over a field that is not algebraically closed. This is the case, for instance, with $\R$ or $\Q_p$. In such situations, the focus shifts to understanding the distribution of the number of solutions when the polynomial system is sampled from a chosen probability measure. This type of question also appears in arithmetic statistics, where often the aim is to determine the proportion of objects, within a specific class, that enjoy a given property. In these contexts, developing efficient methods for sampling from algebraic $p$-adic manifolds is highly desirable.
		
		\begin{figure}[ht]
			\begin{center}
			\scalebox{0.8}{
        			\tikzset{every picture/.style={line width=0.75pt}}
        			
        			\begin{tikzpicture}[x=0.75pt,y=0.75pt,yscale=-0.45,xscale=0.45]
        			
        				\draw [line width=2.25]    (510,10) .. controls (466.19,138.84) and (366.13,181.07) .. (344.28,188.64) .. controls (322.42,196.22) and (310.04,206.34) .. (310,220) .. controls (309.96,233.66) and (317.87,239.93) .. (340,250) .. controls (362.13,260.07) and (462.13,313.07) .. (510,460);
        				
        				\draw  [line width=2.25]  (230,170) .. controls (250.13,186.07) and (253.87,247.93) .. (230,270) .. controls (206.13,292.07) and (149.87,250.93) .. (150,220) .. controls (150.13,189.07) and (209.87,153.93) .. (230,170) -- cycle;
        				
        				\draw  [dash pattern={on 4.5pt off 4.5pt}]  (10,330) -- (600,10);
        				
        				\draw  [dash pattern={on 4.5pt off 4.5pt}]  (140,20) -- (70,440);
        				
        				\draw  [dash pattern={on 4.5pt off 4.5pt}]  (30,60) -- (600,500);
        				
        				\draw  [dash pattern={on 4.5pt off 4.5pt}]  (320,460) -- (230,0);
        				
        				\draw [color={rgb, 250:red, 250; green, 0; blue, 10 }  ,draw opacity=1 ][line width=1.5]    (20,270) -- (630,220);
        				
        				\draw  [fill={rgb, 255:red, 0; green, 7; blue, 255 }  ,fill opacity=1 ] (162.93,257.53) .. controls (162.93,253.37) and (166.31,250) .. (170.47,250) .. controls (174.63,250) and (178,253.37) .. (178,257.53) .. controls (178,261.69) and (174.63,265.07) .. (170.47,265.07) .. controls (166.31,265.07) and (162.93,261.69) .. (162.93,257.53) -- cycle;
        				
        				\draw  [fill={rgb, 255:red, 0; green, 7; blue, 255 }  ,fill opacity=1 ] (325,245) .. controls (325,241.24) and (328.04,238.2) .. (331.8,238.2) .. controls (335.56,238.2) and (338.6,241.24) .. (338.6,245) .. controls (338.6,248.76) and (335.56,251.8) .. (331.8,251.8) .. controls (328.04,251.8) and (325,248.76) .. (325,245) -- cycle;
        				
        				\draw  [color={rgb, 255:red, 36; green, 130; blue, 0 }  ,draw opacity=1 ][fill={rgb, 255:red, 3; green, 149; blue, 0 }  ,fill opacity=1 ] (233.93,252.03) .. controls (233.93,248.15) and (237.08,245) .. (240.97,245) .. controls (244.85,245) and (248,248.15) .. (248,252.03) .. controls (248,255.92) and (244.85,259.07) .. (240.97,259.07) .. controls (237.08,259.07) and (233.93,255.92) .. (233.93,252.03) -- cycle;
        				
        			\end{tikzpicture} }
		
			\end{center}
			\caption{An illustration of the sampling method. The dotted lines are rejected; the red line intersects in 3 points the curve from which we wish to sample. A point is randomly sampled from the three points and the selected point is colored in green.}
		\end{figure}

        For sampling from manifolds over the real or complex numbers, especially in the parameterized case, methods often rely on Markov chain sampling techniques such as the Metropolis algorithm, Gibbs sampler, or hit-and-run algorithm \cite{andersenDiac,saloffCosteDiac,Liu}. While these methods are relatively simple both mathematically and computationally, they only approximate the desired probability measure asymptotically and require a careful analysis of the mixing time of the associated Markov chain \cite{DiacanAl,saloffCosteDiac}. However, $p$-adic fields, being totally disconnected topological spaces, render traditional Markov chain sampling methods ineffective. In contrast, our approach samples exactly from the desired probability density. Furthermore, its geometric nature ensures that it remains effective regardless of the underlying topology.

        \subsection{The setup}

        Let $K$ be a non-archimedean local field. The field $K$ is equipped with a non-archimedean absolute value $|\cdot|$, a valuation map $\val \colon K \twoheadrightarrow \Z \cup \{+\infty\}$, and a unique Haar measure $\mu$ such that $\mu(\Ocal) = 1$, where $\Ocal$ is the unit ball in $K$. We give a brief primer on non-archimedean local fields in \Cref{sec:1}. Readers that are not familiar with the general theory may replace $K$ and $\Ocal$ with $\Q_p$ and $\Z_p$, respectively.
        
        We fix a positive integer $N \geq 2$ and denote by $\A^N := K^{N}$ the $N$-dimensional affine space over $K$. The space $\A^N$ inherits from $K$ the product Haar measure, which we denote by $\mu_{\A^N}$. When there is no risk of confusion, we simply write $\mathrm{d}x$ instead of $\mu_{\A^N}(\mathrm{d}x)$. We endow $\A^N$ with the norm
		\[
		    \norm{x} := \max\limits_{1 \leq i \leq N} |x_i|, \quad \text{for } x = (x_1, \dots, x_N)^{\top} \in \A^N,
		\]
		and the valuation
		\[
		  \val(x) := \min\limits_{1 \leq i \leq N} \val(x_i) = - \log_q(\norm{x}), \quad \text{for } x = (x_1, \dots, x_N)^{\top} \in \A^N,
		\]
        where $q$ is the size of the residue field of $K$ (when $K = \Q_p$, $q$ is equal to $p$). This makes $\A^N$ a metric space with the metric given by $d(x,y) = \norm{x-y}$ for $x,y \in \A^N$. We refer the reader to \Cref{subsec:NormOrthogonality} for more details on norms and valuation.
        
		An affine algebraic variety $X$ in $\A^N$ is the zero set of a system of polynomials $p_1, \dots, p_r$ in $K[x_1, \dots, x_N]$, i.e.
		\[
			X := \big\{  x \in \A^N \colon p_1(x) = \cdots = p_r(x) = 0 \big \}.
		\]
        We refer to smooth and irreducible\footnote{The usual notion of smoothness and irreducibility from algebraic geometry.} varieties over $K$ as \emph{algebraic manifolds}.
		
		Although the notions of dimension and degree of an algebraic $K$-manifold $X$ are the usual notions from algebraic geometry, the notion of \emph{volume} on $X$ is not as standard. Let $X \subset \A^N$ be an affine algebraic $K$-manifold of dimension $n$. For $\epsilon > 0 $ and $x \in \A^N$, let us denote by $B_N(x,\epsilon) = \left\{ y \in \A^N \colon d(x,y) \leq \epsilon \right\}$ the ball of radius $\epsilon$ and center $x$. Similar to the real case, the \emph{volume measure} $\mu_X $ on $X$ is defined as follows: If $V$ is an open set in $X$, then	\begin{equation}\label{eq:defOfAffineMeasure}
			\mu_X(V) \coloneqq  \lim\limits_{\epsilon \to 0} \frac{  \mu_{\A^N} \left(  \bigcup\limits_{ x \in V} B_N(x, \epsilon)  \right) }{ \mu_{\A^{N-n}} \left(  B_{N-n}(0,\epsilon) \right) } = \lim\limits_{r \to \infty}  q^{r (N - n)} \mu_{\A^N} \left(  \bigcup_{x \in V} B_N(x,q^{-r})  \right).
		\end{equation}
		This limit exists (see \cite{Serre81,kulkarni2021integGeo}) and the map $\mu_X$ thus defined\footnote{One could also define a volume measure on $X$ using local charts and differential forms  in the usual way.} is a measure on $X$, where the latter is equipped with its Borel $\sigma$-algebra.
		
		\begin{remark}
		  When $X \subset \P^{N-1}$ is a projective manifold, one can still define a volume measure $\mu_X$ in the same manner by replacing the measure $\mu_{\A^N}$ in (\ref{eq:defOfAffineMeasure}) with its normalized push-forward to the projective space $\P^{N-1}$, and by defining balls in $\P^{N-1}$ using the \emph{Fubini–Study metric}. Here we focus on the affine case and delay our treatment of projective manifolds until \Cref{sec:Sampling_ProjectiveManifolds}.
		\end{remark}

		Given a function $f \colon X \to \R$, integrable with respect to the measure $\mu_X$, we wish~to:
		\begin{enumerate}
		    \item Estimate the integral $\int_{X} f(x) \mu_X(\mathrm{d}x)$.
		    \item Sample a random variable $\bm{\xi} \in X$ with the probability density $f(x) / \int_{X} f(x) \mu_X(\mathrm{d}x)$, when $f$ is non-negative and $\int_{X} f(x) \mu_X(\mathrm{d}x) > 0$.
		\end{enumerate}
		
        \medskip
		
        Our sampling approach is based on intersecting the manifold $X$ with affine linear subspaces of complementary dimension. For a matrix $A \in K^{n \times N}$ and $b \in K^n$, we define the affine linear space $\Lcal_{A,b}$ as follows:  
        \[
        \Lcal_{A,b} \coloneqq \{ x \in \A^N \colon Ax = b \}.
        \]  
        This space is generically\footnote{The set of pairs $(A, b) \in K^{n \times N} \times K^n$ for which $\Lcal_{A,b}$ has dimension $N - n$ is non-empty and Zariski-open.} of dimension $N-n$, and the intersection $\Lcal_{A,b} \cap X$ is generically\footnote{The set of pairs $(A, b) \in K^{n \times N} \times K^n$ for which the intersection is finite is non-empty and Zariski-open.} finite with its size bounded above by the degree of $X$ (see \Cref{subsec:intersection}). 
        
        In essence, sampling from $X$ can be reduced to two steps: sampling a random affine space $\Lcal_{A, b}$, and then sampling a random point from the finite set $\Lcal_{A, b} \cap X$. However, given a target probability density $f$ on $X$, neither of these steps is entirely straightforward. For this reason, we introduce a \textit{weight} function 
        \[
            w_X \colon X \to \R_{ > 0}.
        \]
        The precise definition of $w_X$ is deferred to \Cref{subsec:weighFunc}.

        \subsection{Main results}
        
        Given a real valued function $f$ on $X$, we define the following function:
		\[
    		\overline{f} \colon K^{n \times N}\times K^n \xrightarrow[]{} \R,  \quad  (A,b) \mapsto \sum\limits_{x \in X \cap \Lcal_{A,b}} w_X(x) f(x).
		\]
		When $X \cap \Lcal_{A,b}$ is empty or infinite, we set $\overline{f}(A,b) = 0$. Our first result deals with integrating a real-valued integrable function $f$ on an algebraic manifold $X$. Namely, we show that the integral can be expressed as the expectation of a real-valued random variable that we can sample.
		
		\begin{theorem}
		\label{thm:integral_EP}
		Let $X \subset \A^{N}$ be an $n$-dimensional affine algebraic manifold defined over $K$. Let $(\bm{A},\bm{b})$ be a random variable in $K^{n \times N} \times K^n$ with independent entries uniformly distributed in $\Ocal$, i.e., $(\bm{A},\bm{b})$ have distribution $1_{A \in \Ocal^{n \times N}, \ b \in \Ocal^n}  \mathrm{d}A \mathrm{d}b$. Then we have:
		\[
		\int_{X} f(x) \mu_{X}(\mathrm{d}x) = \E\left[\overline{f}(\bm{A},\bm{b})\right].
		\]
		\end{theorem}
		
		With this theorem in hand, one can evaluate integrals, up to a certain confidence interval, using Monte-Carlo methods. We discuss this in more detail in \Cref{sec:Sampling_practically}. 
        
        \medskip	
        
		Our second result deals with sampling a random point $\bm{\xi}$ from  a manifold $X$ with a prescribed probability density $f$ with respect to the natural volume measure $\mu_X$ on $X$:
		
		\begin{theorem} \label{thm:xi_density}
			Let $X \subset \A^{N}$ be an $n$-dimensional affine algebraic manifold defined over $K$. Let $f \colon X \to \R_{\geq 0}$ be a probability density with respect to $\mu_X$. Let $(\bm{\widetilde{A}}, \bm{\widetilde{b}})$ be the random variable in $K^{n \times N} \times K^n$ with distribution
			\[
			        \overline{f}(A,b) \ 1_{A \in \Ocal^{n \times N}, \ b \in \Ocal^n} \ \mathrm{d}A \mathrm{d}b.
			\]
			Let $\bm{\xi}$ be the random variable obtained by intersecting $X$ with the random space $\Lcal_{\bm{\widetilde{A}}, \bm{\widetilde{b}}}$ and choosing a point $x$ in the finite set $X \cap \Lcal_{\bm{\widetilde{A}}, \bm{\widetilde{b}}}$ with probability 
			\[
			    \frac{w_X(x)f(x)}{\overline{f}(\bm{\widetilde{A}}, \bm{\widetilde{b}})}.
			\]
			Then $\bm{\xi}$ has density $f$ with respect to the volume measure $\mu_X$ on $X$. 
		\end{theorem}
		 
		We give similar results for projective manifolds, namely Theorems \ref{thm:integral_EPprojective} and \ref{thm:xi_densityProjective}, in \Cref{sec:Sampling_ProjectiveManifolds}. We provide a SageMath \cite{sagemath} implementation of the sampling method we describe in this article (in some particular cases) in the repository \cite{MathRepo}.
		
		\begin{remark}
		    Although most of our results are stated for algebraic manifolds, there is no issue working with an irreducible variety $X$ (affine or projective) with potential singularities. This is because the singular locus $X^{\mathrm{sing}}$ is lower dimensional in $X$ and we can work with the algebraic manifold $X \setminus X^{\mathrm{sing}}$. Our sampling method will then produce a point in $X$ that is smooth with probability $1$.
		\end{remark}

        \subsection{Notation}
		
		For the reader's convenience, we summarize some of the notation we use throughout this paper in the following table.
        
        \medskip
		
				\begin{center}
						\begin{tabular}{rcl}
						    $K$         		   		   			 &--& 	A non-archimedean local field.\\
							$\val(\cdot)$  	    					 &--&    The valuation on $K$ and also, depending on context, the valuation on $K^N$.\\
							$\Ocal$									 &--&	The valuation ring of $K$.\\
							$\varpi$                                 &--&    A uniformizer of $K$; that is $\val(\varpi) = 1$.\\
							$k$                                      &--&    The residue field $\Ocal/\varpi \Ocal$.\\
                            $q$                                      &--&    The size of the residue field $k$.\\
                            $p$                                      &--&    The characteristic of the residue field $k$.\\
							$|\cdot|$   		   		   			 &--&	The absolute value on $K$ given by $|x| := q^{-\val(x)}$ for any $x \in K$.\\
							$\norm{\cdot}$  	    			     &--&    The standard norm on the vector space $K^{N}$, associated with the standard \\
																     &  &    lattice $\Ocal^N$, see \Cref{subsec:NormOrthogonality}.\\
							$X$           		   					 &--& 	An algebraic $K$-manifold (either affine or projective). \\
							$w_X$          		  	      			 &--&    The weight function on $X$ when $X$ is affine. \\
							$\mu_X$          		  	      		 &--&    The volume measure on $X$. \\
							$1_S$                                    &--&    The indicator of a set $S$; that is $1_{x \in S} = 1_S(x) = 1$ if $x \in S$ and $0$ otherwise.
						\end{tabular}
				\end{center}
                \medskip

                Finally, for clarity, we write all random variables in bold font.

		\section{Background and preliminaries}\label{Sec:2}
        
		In this section, we expand on the concepts used in the statements of Theorems \ref{thm:integral_EP} and \ref{thm:xi_density}, and gather some of the tools that will be needed throughout the remainder of the article. The following subsection is essentially a brief summary of standard material that can be found in a number of sources, for example \cite{SerreLocalField}.

        \subsection{A primer on non-archimedean local fields}\label{sec:1}
        
        We begin with the prototypical example: the field of \emph{$p$-adic numbers}~$\Q_p$. Let $p$ be a fixed prime number. For any nonzero $x \in \Q$, there exists a unique integer $\val_p(x) \in \Z$ such that  
        \[
        x = p^{\val_p(x)} \cdot \frac{a}{b}, \quad \text{where } a, b \in \Z \setminus p\Z.
        \]  
        By setting $\val_p(0) = +\infty$, we obtain the \emph{$p$-adic valuation}, $\val_p : \Q \to \Z \cup \{+\infty\}$. The real-valued map  
        \[
        |\cdot|_p : \Q \to \mathbb{R}_{\geq 0}, \qquad |x|_p := p^{-\val_p(x)},
        \]  
        is called the \emph{$p$-adic absolute value}. This absolute value is non-archimedean, meaning it satisfies the following properties for all $x, y \in \Q$:  
        
        \begin{enumerate}[wide=80pt, leftmargin=45pt]
            \item $|x|_p = 0$ if and only if $x = 0$,
            \item $|xy|_p = |x|_p \cdot |y|_p$, \hfill (Multiplicativity)
            \item $|x+y|_p \leq \max(|x|_p, |y|_p)$. \hfill (Ultrametric inequality)
        \end{enumerate}

        \medskip
        \noindent The map $(x,y) \mapsto |x-y|_p$ defines an ultrametric on $\Q$ and the completion of $\Q$ with respect to this metric is the field of $p$-adic numbers $\Q_p$ and $|\cdot|_p$ naturally extends to $\Q_p$. The field $\Q_p$ then has the structure of a topological field.
    
       \begin{example}
            Let $p = 7$ and consider the sequence $(x_n := 1 + 7 + \dots + 7^n)_{n \geq 0}$. This sequence is Cauchy in $\Q$ with respect to the $7$-adic absolute value $|\cdot|_7$ because 
            \[
                |x_n - x_m|_7 = 7^{-(m+1)} \quad \text{for all } n \geq m \geq 0.
            \]
            Thus, the sequence $(x_n)$ converges to a $7$-adic number $x = 1 + 7 + 7^2 + \dots \in \Q_7$. In fact, $x = -6$, since $6 + x \in \Z_7$ is divisible in $\Z_7$ by $7^n$ for all $n \geq 0$. On the other hand, the sequence $(y_n := 1 + 7^{-1} + \dots + 7^{-n})_{n \geq 0}$ does not converge in $\Q_7$, because
            \[
                |y_{n+1} - y_n|_7 = 7^{n+1} \quad \text{for all } n \geq 1.
            \]            
        \end{example}
        
          The properties of $|\cdot|_p$ ensure that the unit ball  
        \[
        \Z_p := \{ x \in \Q_p \colon |x|_p \leq 1 \}
        \]  
        is compact and has the structure of a ring. This ring is called the \emph{valuation ring} of $\Q_p$. Alternatively, $\Z_p$ can be constructed as a topological ring via the projective limit  
        \[
        \Z_p := \varprojlim \Z/p^n \Z,
        \]  
        where each factor $\Z/p^n\Z$ is equipped with the discrete topology. Here, the limit is taken with respect to the natural morphisms $\Z/p^{n}\Z \to \Z/p^{m}\Z$ for $n \geq m$. The field $\Q_p$ can then be recovered as the field of fractions of $\Z_p$. This construction implies that any nonzero element $x \in \Q_p$ can be expressed as a Laurent series in $p$:  
        \[
        x = \sum_{n = \val_p(x)}^{\infty} a_n p^n, \quad a_n \in \{0, 1, 2, \dots, p-1\} \text{ with } a_{\val_p(x)} \neq 0.
        \]  
        
    \begin{figure}
        \centering
        \scalebox{0.5}{
                \tikzset{every picture/.style={line width=0.75pt}}       
                
                \begin{tikzpicture}[x=0.75pt,y=0.75pt,yscale=-1,xscale=1]
                
                \draw  [fill={rgb, 255:red, 245; green, 166; blue, 35 }  ,fill opacity=0.3 ] (92,237.75) .. controls (92,107.55) and (197.55,2) .. (327.75,2) .. controls (457.95,2) and (563.5,107.55) .. (563.5,237.75) .. controls (563.5,367.95) and (457.95,473.5) .. (327.75,473.5) .. controls (197.55,473.5) and (92,367.95) .. (92,237.75) -- cycle ;
                \draw  [fill={rgb, 255:red, 0; green, 47; blue, 255 }  ,fill opacity=0.3 ] (115.3,317.3) .. controls (115.3,265.27) and (160,223.09) .. (215.15,223.09) .. controls (270.3,223.09) and (315,265.27) .. (315,317.3) .. controls (315,369.32) and (270.3,411.5) .. (215.15,411.5) .. controls (160,411.5) and (115.3,369.32) .. (115.3,317.3) -- cycle ;
                \draw  [fill={rgb, 255:red, 208; green, 2; blue, 27 }  ,fill opacity=0.6 ] (123.42,343.55) .. controls (123.42,320.24) and (143.45,301.34) .. (168.16,301.34) .. controls (192.87,301.34) and (212.9,320.24) .. (212.9,343.55) .. controls (212.9,366.86) and (192.87,385.75) .. (168.16,385.75) .. controls (143.45,385.75) and (123.42,366.86) .. (123.42,343.55) -- cycle ;
                \draw  [fill={rgb, 255:red, 208; green, 2; blue, 27 }  ,fill opacity=0.6 ] (169,264.41) .. controls (169,241.1) and (189.03,222.2) .. (213.74,222.2) .. controls (238.44,222.2) and (258.47,241.1) .. (258.47,264.41) .. controls (258.47,287.72) and (238.44,306.62) .. (213.74,306.62) .. controls (189.03,306.62) and (169,287.72) .. (169,264.41) -- cycle ;
                \draw  [fill={rgb, 255:red, 208; green, 2; blue, 27 }  ,fill opacity=0.6 ] (219.13,341.29) .. controls (219.13,317.98) and (239.16,299.09) .. (263.87,299.09) .. controls (288.58,299.09) and (308.61,317.98) .. (308.61,341.29) .. controls (308.61,364.6) and (288.58,383.5) .. (263.87,383.5) .. controls (239.16,383.5) and (219.13,364.6) .. (219.13,341.29) -- cycle ;
                \draw  [fill={rgb, 255:red, 0; green, 47; blue, 255 }  ,fill opacity=0.3 ] (339.3,316.3) .. controls (339.3,264.27) and (384,222.09) .. (439.15,222.09) .. controls (494.3,222.09) and (539,264.27) .. (539,316.3) .. controls (539,368.32) and (494.3,410.5) .. (439.15,410.5) .. controls (384,410.5) and (339.3,368.32) .. (339.3,316.3) -- cycle ;
                \draw  [fill={rgb, 255:red, 208; green, 2; blue, 27 }  ,fill opacity=0.6 ] (346.42,342.55) .. controls (346.42,319.24) and (366.45,300.34) .. (391.16,300.34) .. controls (415.87,300.34) and (435.9,319.24) .. (435.9,342.55) .. controls (435.9,365.86) and (415.87,384.75) .. (391.16,384.75) .. controls (366.45,384.75) and (346.42,365.86) .. (346.42,342.55) -- cycle ;
                \draw  [fill={rgb, 255:red, 208; green, 2; blue, 27 }  ,fill opacity=0.6 ] (392,263.41) .. controls (392,240.1) and (412.03,221.2) .. (436.74,221.2) .. controls (461.44,221.2) and (481.47,240.1) .. (481.47,263.41) .. controls (481.47,286.72) and (461.44,305.62) .. (436.74,305.62) .. controls (412.03,305.62) and (392,286.72) .. (392,263.41) -- cycle ;
                \draw  [fill={rgb, 255:red, 208; green, 2; blue, 27 }  ,fill opacity=0.6 ] (442.13,340.29) .. controls (442.13,316.98) and (462.16,298.09) .. (486.87,298.09) .. controls (511.58,298.09) and (531.61,316.98) .. (531.61,340.29) .. controls (531.61,363.6) and (511.58,382.5) .. (486.87,382.5) .. controls (462.16,382.5) and (442.13,363.6) .. (442.13,340.29) -- cycle ;
                \draw  [fill={rgb, 255:red, 0; green, 47; blue, 255 }  ,fill opacity=0.3 ] (225.3,98.3) .. controls (225.3,46.27) and (270,4.09) .. (325.15,4.09) .. controls (380.3,4.09) and (425,46.27) .. (425,98.3) .. controls (425,150.32) and (380.3,192.5) .. (325.15,192.5) .. controls (270,192.5) and (225.3,150.32) .. (225.3,98.3) -- cycle ;
                \draw  [fill={rgb, 255:red, 208; green, 2; blue, 27 }  ,fill opacity=0.6 ] (233.42,124.55) .. controls (233.42,101.24) and (253.45,82.34) .. (278.16,82.34) .. controls (302.87,82.34) and (322.9,101.24) .. (322.9,124.55) .. controls (322.9,147.86) and (302.87,166.75) .. (278.16,166.75) .. controls (253.45,166.75) and (233.42,147.86) .. (233.42,124.55) -- cycle ;
                \draw  [fill={rgb, 255:red, 208; green, 2; blue, 27 }  ,fill opacity=0.6 ] (280,45.41) .. controls (280,22.1) and (300.03,3.2) .. (324.74,3.2) .. controls (349.44,3.2) and (369.47,22.1) .. (369.47,45.41) .. controls (369.47,68.72) and (349.44,87.62) .. (324.74,87.62) .. controls (300.03,87.62) and (280,68.72) .. (280,45.41) -- cycle ;
                \draw  [fill={rgb, 255:red, 208; green, 2; blue, 27 }  ,fill opacity=0.6 ] (329.13,122.29) .. controls (329.13,98.98) and (349.16,80.09) .. (373.87,80.09) .. controls (398.58,80.09) and (418.61,98.98) .. (418.61,122.29) .. controls (418.61,145.6) and (398.58,164.5) .. (373.87,164.5) .. controls (349.16,164.5) and (329.13,145.6) .. (329.13,122.29) -- cycle ;
                
                \draw (307,29.4) node [anchor=north west][inner sep=0.75pt] [xscale=1.4, yscale=1.4]    {$3^{2}  \Z_{3}$};
                \draw (238,114.4) node [anchor=north west][inner sep=0.75pt] [xscale=1.4, yscale=1.4]   {$3+3^{2}  \Z_{3}$};
                \draw (336,109.4) node [anchor=north west][inner sep=0.75pt] [xscale=1.4, yscale=1.4]   {$6+3^{2}  \Z_{3}$};
                
                \draw (175,248.4) node [anchor=north west][inner sep=0.75pt] [xscale=1.4, yscale=1.4]   {$1+3^{2}  \Z_{3}$};
                \draw (128,328) node [anchor=north west][inner sep=0.75pt] [xscale=1.4, yscale=1.4]   {$4+3^{2}  \Z_{3}$};
                \draw (226,328) node [anchor=north west][inner sep=0.75pt] [xscale=1.4, yscale=1.4]   {$7+3^{2}  \Z_{3}$};

                \draw (395,250) node [anchor=north west][inner sep=0.75pt][xscale=1.4, yscale=1.4]    {$2+3^{2}  \Z_{3}$};
                \draw (348,332.4) node [anchor=north west][inner sep=0.75pt] [xscale=1.4, yscale=1.4]   {$5+3^{2}  \Z_{3}$};
                \draw (448,331.4) node [anchor=north west][inner sep=0.75pt] [xscale=1.4, yscale=1.4]   {$8+3^{2}  \Z_{3}$};

                \draw (180,53.4) node [anchor=north west][inner sep=0.75pt] [xscale=1.6, yscale=1.6]   {$3  \Z_{3}$};
                \draw (90,210) node [anchor=north west][inner sep=0.75pt] [xscale=1.6, yscale=1.6]   {$1+3  \Z_{3}$};
                \draw (470,205) node [anchor=north west][inner sep=0.75pt] [xscale=1.6, yscale=1.6]   {$2+3  \Z_{3}$};

                \draw (321,476.4) node [anchor=north west][inner sep=0.75pt] [xscale=2, yscale=2]   {$ \Z_{3}$};
                \end{tikzpicture}
        }
            
        \caption{An illustration of the $3$-adic ring $\Z_3$. This picture shows all of the balls in $\Z_3$ of radius $r$ such that $3^{-2} \leq r \leq 1$. These are of the form $a + 3^k \Z_3$ where $ 0 \leq k \leq 2$ and $a \in \{ 0, 1, \dots, 3^{2} - 1 \}$. Because of the ultrametric nature of the $3$-adic absolute value, any two balls of the same radius are either equal or disjoint. This is a feature of the $p$-adic topology.}
        \label{fig:enter-label}
    \end{figure}

        The ring $\Z_p$ is a local ring with a unique maximal ideal $p \Z_p$. Any generator of this ideal, such as $p$, is referred to as a \emph{uniformizer} of $\Q_p$. The residue field of $\Q_p$ is $\F_p := \Z_p / p \Z_p = \Z / p \Z$, which has characteristic $p$. In contrast, the field $\Q_p$ itself has characteristic $0$. This interplay of characteristics makes $\Q_p$ a local field of \emph{mixed characteristic}.
        
        Since the group $p^{\Z}$ is discrete in $\mathbb{R}_{>0}$, the ball $\Z_p$ is both open and closed (clopen). The topology of $\Q_p$ has a basis of clopen sets consisting of translates of the form $a + p^k \Z_p$, where $a \in \Q_p$ and $k \in \Z$. As a result, $\Q_p$ is a locally compact, totally disconnected, non-discrete topological field.
        
        The field $\Q_p$ is endowed with its \emph{Borel $\sigma$-algebra} $\mathcal{B}(\Q_p)$, which is the $\sigma$-algebra generated by the open sets of $\Q_p$. A \emph{Haar measure} on $\Q_p$ is a measure $\mu \colon \mathcal{B}(\Q_p) \to \R_{\geq 0}$ satisfying the following properties for any Borel set $A \subset \Q_p$ and $x \in \Q_p$:
        \vspace{2mm}
        \begin{enumerate}[wide=100pt, leftmargin=45pt]
            \item $\mu(x + A) = \mu(A)$ \hfill (translation invariance),
            \item $\mu(x \cdot A) = |x|_p \cdot \mu(A)$ \hfill (scaling property),
        \end{enumerate}
        where 
        \[
        x + A = \{ x + a \colon a \in A \} \quad \text{and} \quad x \cdot A = \{ x a \colon a \in A \}.
        \]
        Since $\Q_p$ is locally compact, such a measure is unique up to scalar multiplication. We denote by $\mu$ the unique Haar measure on $\Q_p$ that satisfies $\mu(\Z_p) = 1$. In particular, the measure $\mu$ assigns to each ball of the form $p^k \Z_p$ the value
        \[
                \mu\left(p^k \Z_p\right) = p^{-k}.
        \]
        As $k$ tends to $+ \infty$, the ball $p^k \Z_p$ shrinks and its measure converges to zero. The measure $\mu$ is a probability measure on $\Z_p$, and a random element $\bm{\xi} \in \Z_p$ with distribution $\mu$ can be generated as
        \[
        \bm{\xi} = \sum_{n = 0}^{\infty} \bm{a_n} p^n,
        \]
        where $(\bm{a_n})$ is a sequence of independent and identically distributed (i.i.d.) random variables, each uniformly distributed in $\{0, \dots, p-1\}$.
    
        \medskip
        
       Another important example of a non-archimedean local field is the field $\F_p((t))$ of formal Laurent series in the variable $t$ with coefficients in the finite field $\F_p$. Any nonzero element $f \in \F_p((t))$ can be uniquely written as a formal series
        \[
        f = \sum_{n = m}^{+\infty} a_n t^n, \quad m \in \Z, \quad a_n \in \F_p \text{ for all } n \geq m, \quad \text{and} \quad a_m \neq 0.
        \]
        The $t$-adic valuation of $f$ is defined as $\val_t(f) := m$, and the absolute value is $|f|_t := p^{-\val_t(f)}$. The unit ball (or valuation ring) of $\F_p((t))$ is the ring $\F_p[[t]]$ of formal power series in $t$, which is a local ring with maximal ideal $t \F_p[[t]]$, generated by the uniformizer $t$. Its residue field is $\F_p$, i.e., $\F_p[[t]] / t \F_p[[t]] = \F_p$. Both $\F_p((t))$ and its residue field $\F_p$ have characteristic $p$, making $\F_p((t))$ a local field of \emph{equal characteristic}.
        
        Similar to $\Q_p$, there exists a unique Haar measure on $\F_p((t))$ that assigns measure 1 to the unit ball $\F_p[[t]]$ and a random element $\bm{\xi}$ in $\F_p[[t]]$ with distribution $\mu$ can be generated as
        \[
            \bm{\xi} = \sum_{n = 0}^{\infty} \bm{a_n} t^n,
        \]
        where $(\bm{a_n})_{n \geq 0}$ is an independent and identically distributed (i.i.d) sequence of random variables, each uniformly distributed in $\F_p$.
        
        \medskip
        At this stage, readers unfamiliar with general non-archimedean local fields may choose to skip the remainder of this subsection or consult it only for notation. For the rest of the article, it is sufficient to work with either of the local fields $K = \Q_p$ or $K = \F_p((t))$.

        In general, a non-archimedean local field is a locally compact, non-discrete, and totally disconnected topological field. For the remainder of this text, we fix such a field $K$. It is a well-known fact that $K$ is isomorphic to a finite extension of $\Q_p$ or $\F_p((t))$ for some prime $p$. For further details and background, we refer the reader to \cite[Chapter II]{SerreLocalField}. Specifically, the following properties hold:
        
        \smallskip
        
        \begin{enumerate}[wide=40pt, leftmargin=54pt]
            \item There exists a unique normalized discrete valuation $\val \colon K \twoheadrightarrow \Z \cup \{+\infty\}$ satisfying the following properties for all $x, y \in K$:
            \begin{enumerate}[wide=100pt]
                \item $\val(x) = +\infty$ if and only if $x = 0$,
                \item $\val(xy) = \val(x) + \val(y)$
                \item $\val(x+y) \geq \min(\val(x), \val(y))$.
            \end{enumerate}
            \vspace{2mm}
            
            \item The valuation ring $\Ocal := \{x \in K \colon \val(x) \geq 0\}$ is a local ring with a unique maximal ideal $\varpi \Ocal$, where $\varpi \in K$ is any uniformizer, i.e., an element satisfying $\val(\varpi) = 1$.
        
            \item The residue field $k := \Ocal / \varpi \Ocal$ is finite with size $q := |k|$, which is a power of the residue characteristic $p := \Char(k)$. 
        
            \item The absolute value $x \mapsto |x| := q^{-\val(x)}$ is ultrametric, and $K$ is complete with respect to $|\cdot|$.
        \end{enumerate}
        
        \begin{remark}
        For each nonzero element $\alpha \in k = \Ocal / \varpi \Ocal$, choose a lift $\widetilde{\alpha} \in \Ocal$ and define the set of representatives $S := \{\widetilde{\alpha} \colon \alpha \in k^\times\} \cup \{0\}$. Similar to $\Q_p$ and $\F_p((t))$, any nonzero element $x \in K$ can be uniquely expressed as
        \[
        x = \sum_{n = \val(x)}^{+\infty} a_n \varpi^n, \quad a_n \in S \text{ for all } n \geq \val(x), \quad \text{and} \quad a_{\val(x)} \neq 0.
        \]
        \end{remark}
        
        Since $K$ is locally compact, it admits a unique \emph{Haar measure} $\mu \colon \Bcal(K) \to \R_{\geq 0}$ normalized such that $\mu(\Ocal) = 1$, where $\Bcal(K)$ denotes the collection of Borel subsets of $K$. A random element $\bm{\xi}$ in $\Ocal$ with distribution $\mu$ can be generated as
        \[
        \bm{\xi} = \sum_{n=0}^\infty \bm{a}_n \varpi^n,
        \]
        where $(\bm{a}_n)_{n \geq 0}$ is a sequence of i.i.d. random variables, each uniformly distributed over the finite set of representatives $S$.

        \subsection{Norms, Lattices, and Gaussian Measures} \label{subsec:NormOrthogonality}
        
        Let \( E \) be a vector space over \( K \) with \( N = \dim_K(E) \). A \emph{lattice} in \( E \) is an \( \mathcal{O} \)-submodule \( \Lambda \subseteq E \) of rank \( N \). More explicitly, a lattice \( \Lambda \) is a module generated by a basis \( (e_1, \dots, e_N) \) of \( E \) over \( \mathcal{O} \), i.e.
        \[
        \Lambda = \mathcal{O} e_1 \oplus \cdots \oplus \mathcal{O} e_N.
        \]
        Each lattice \( \Lambda \) induces a valuation map \( \val_\Lambda \) on \( E \), defined by
        \[
        \val_\Lambda(x) \coloneqq \sup \left\{ m \in \mathbb{Z} \colon \varpi^{-m} x \in \Lambda \right\}, \quad x \in E.
        \]
        This valuation satisfies the following properties:
        \begin{align*}
        (1) & \quad \val_\Lambda(x) = + \infty \iff x = 0, \quad \text{for all } x \in E, \\
        (2) & \quad \val_\Lambda(\alpha x) = \val(\alpha) + \val_\Lambda(x), \quad \text{for all } \alpha \in K, x \in E, \\
        (3) & \quad \val_\Lambda(x + y) \geq \min(\val_\Lambda(x), \val_\Lambda(y)),  \quad \text{for all } x,y \in E.
        \end{align*}
        Using this valuation, we can define an ultrametric norm \( \norm{\cdot}_\Lambda \) on the vector space \( E \) by
        \[
        \norm{x}_\Lambda \coloneqq q^{-\val_\Lambda(x)}, \quad x \in E.
        \]
        Note that the lattice \( \Lambda \) corresponds to the unit ball with respect to this norm, i.e.
        \[
        \Lambda = \{ x \in E \colon \norm{x}_\Lambda \leq 1 \}.
        \]
        Thus, there is a bijection between norms on \( E \) and lattices in \( E \).
        
        \medskip
         
         Given a lattice $\Lambda$ in $E$ (or a norm on $E$), there is a natural notion of orthogonality on $(E, \Lambda, \norm{\cdot}_\Lambda)$ for which a family of vectors $(e_i)_{i \in I}$ is \emph{orthogonal} if
         \[
            \norm{ \sum_{j \in J}  \alpha_j e_j }_\Lambda = \max_{j \in J} |\alpha_j| \norm{e_j}_\Lambda, \quad \text{ for any finite set } J \subset I \text{ and any } \alpha_j \in K.
         \]
         Fixing an orthonormal basis of $e_1, \dots, e_N$ of $E$, the group of linear automorphisms of $E$ preserving the norm $\norm{\cdot}_\Lambda$, that is 
         \[
            G_{\Lambda} \coloneqq \left\{ g \in \GL(E) \colon \norm{g(x)}_{\Lambda}  = \norm{x}_\Lambda \right\}
         \]
         is isomorphic the group of matrices $\GL(N,\Ocal)$.

         When $E = K^N$ we call $\Ocal^N$ the \emph{standard lattice in $K^N$}, i.e. the lattice generated by the standard basis of $K^N$ over $\Ocal$, and we denote by $\norm{\cdot}$ its associated norm. Explicitly, this norm and its corresponding valuation are given by
         \begin{equation} \label{eq:norm_on_KN}
            \val(x) \coloneqq \min\limits_{1 \leq i \leq N} \val(x_i) \text{ and } \norm{x} \coloneqq \max\limits_{1 \leq i \leq N} |x_i| ,\quad \text{ for } x = (x_1, \dots, x_N)^\top \in K^N.
         \end{equation}
         From now on, unless mentioned explicitly, the norm on $K^N$ is the norm given by \Cref{eq:norm_on_KN}. 
         
         The group $\GL(N,\Ocal)$ is the group of matrices with orthonormal rows and columns. This is also the group of linear isometries of $E = K^N$ and the Haar measure $\mu_{\A^N}$ is invariant under its natural action $\A^N$, i.e.,
          \[
                \mu_{\A^N} (g \cdot A ) = \mu_{\A^N}(A), \quad \text{for any } g \in \GL(N,\Ocal) \text{ and Borel set } A \subset K^N.
          \]
          The uniform probability measure on the standard lattice $\Ocal^N$ (with respect to the Haar measure) is called the \emph{standard Gaussian distribution} on $E$.
        
          For readers who are not familiar with non-archimedean orthogonality and Gaussian measures, we recommend consulting \cite{evans2001local} as well as \cite{maazouz2019statistics,maazouz2021gaussian} for recent overviews. For a more comprehensive treatment of non-archimedean functional analysis, we refer to \cite[Chapter 5]{vanRooij}.

           \newpage

        \subsection{The weight function of an affine algebraic manifold}\label{subsec:weighFunc}

        We now define the \emph{weight function} we mentioned in the introduction. Let $X \subset \A^N$ be an affine algebraic manifold defined over $K$ and let $n := \dim(X)$. However, before we do so, we begin with a couple of definitions.
		
		\begin{definition}\label{def:absoluteDet}
		Let $a,b \geq 1$ be two positive integers, and $M \in K^{a \times b}$ a matrix. Let us write the Smith normal form of $M$ as $M = U D V$, where $U \in \GL(a,\Ocal); V \in \GL(b,\Ocal)$; and $D = \diag(\varpi^{v_1}, \dots, \varpi^{v_{\min(a,b)}}) \in K^{a \times b}$ with $v_1 \geq \cdots \geq v_{\min(a,b)} \in \Z \cup \{ \infty \}$. We then define the \emph{absolute determinant} of $M$ as follows:
		\[
		    \N(M) \coloneqq |\varpi^{v_1} \cdots \varpi^{v_{\min(a,b)}}| = q^{- v_1 - \cdots - v_{\min(a,b)} }.
		\]
        The diagonal elements $\varpi^{-v_1}, \dots, \varpi^{-v_{\min(a,b)}}$ are called the \emph{singular values} of the matrix $M$.
  
		If $E,F$ are $K$-vector spaces of respective dimensions $b,a$ and $\varphi \colon E \to F$ is $K$-linear, we define
		\[
		\N(\varphi) = \N(A),
		\]
		where $A \in K^{a \times b}$ is a matrix representing $\varphi$ in orthonormal bases\footnote{In the sense of \Cref{subsec:NormOrthogonality}.} of $E$ and $F$.
		\end{definition}

		\begin{definition}\label{def:NrDefAffine}
		Let $X \subset \A^N$ be an affine algebraic manifold of dimension $n$, and $x$ be a point on $X$. Let $U \in \GL(N,\Ocal)$ such that $U x = (0, \dots,0, \varpi^{\val(x)})^\top$, and let $W \in \Ocal^{N \times n}$ be a matrix whose columns form an orthonormal basis of the tangent space $T_xX$. Finally, let us set $S_x = \diag(1, \dots, 1, \varpi^{\max(0, - \val(x))}) \in K^{N \times N}$. Then we define 
		\[
			\Nr(X,x) \coloneqq \N(S_x U W).
		\]
		\end{definition}
        
		It is not so clear apriori that this definition does not depend on the choice of $W$ and $U$. The following lemma proves that quantity $\Nr(X,x)$ is independent of this choice.

		\begin{lemma}\label{lem:indepOfChoice}
		    Let $X \subset \A^N$ be an affine algebraic $K$-manifold of dimension $n$. The definition of $\Nr(X,x)$ in \Cref{def:NrDefAffine} does not depend on the choice of $U$ and $W$.
		\end{lemma}
		\begin{proof}
		Set $S_x = \diag(1, \dots, 1, \varpi^{\max(0, -\val(x))})$ and let $U_1, U_2 \in \GL(N,\Ocal)$ and $W_1, W_2 \in \Ocal^{N \times n}$ be such that
		\begin{enumerate}[label=(\arabic*)]
		    \item $U_1 x =  U_2 x = (0, \dots , 0, \varpi^{\val(x)})^\top$, \label{cond:1}
		    \item columns of $W_1,W_2$ are two orthonormal bases of the tangent space $T_x X$.
		\end{enumerate}
		Then there exists $V \in \GL(n,\Ocal)$ such that $W_2 = W_1V$. Let $A = U_2 U_1^{-1}$ and $B = S_x A S_{x}^{-1} $ so that  we have $B (S_x U_1 W_1) V = S_x U_2 W_2$. We claim that $B \in \GL(N,\Ocal)$ is an orthogonal matrix. To see why, notice that, thanks to condition \ref{cond:1}, the matrix $A$ is of the form
		\[
		A = \begin{pmatrix}
              A' & \rvline & \begin{matrix} 0\\ \vdots \\ 0 \end{matrix}  \\ 
                   \hline
             \begin{matrix} z_1 & \cdots & z_{N-1} \end{matrix}     & \rvline & 1
        \end{pmatrix}
		\]
		 where $A' \in \GL(N-1, \Ocal)$ and $z_1, \dots, z_{N-1} \in \Ocal$. We then deduce that $B$ is of the form
		 \[
		 B = \begin{pmatrix}
              A' & \rvline & \begin{matrix} 0\\ \vdots \\ 0 \end{matrix}  \\ 
                   \hline
             \begin{matrix} \alpha z_1 & \cdots & \alpha z_{N-1} \end{matrix}     & \rvline & 1
        \end{pmatrix}
		 \]
		 where $\alpha = \varpi^{\max(0, -\val(x))} \in \Ocal$. So we deduce that $B \in \GL(N,\Ocal)$. Now, since $V$ and $B$ are both orthogonal, from \Cref{def:absoluteDet} we can see that
		 \[
		    \N(S_x U_1W_1) = \N(S_x U_2 W_2) = \N(B S_x U_1 W_1 V),
		 \]
		 which finishes the proof.
		\end{proof}    

        \begin{remark}
            The quantity $\Nr(X,x)$ can be interpreted as a ``measure" of how far the point $x$ is from being an $\Ocal$-point of $X$.     
        \end{remark}
        
		Now we are ready to define the weight function $w_X$ on $X$.
		
		\begin{definition}\label{def:affineWeight}
		Let $X \subset \A^{N}$ be an affine algebraic manifold of dimension $n$ over $K$. For a point $x \in X$, we define the \emph{weight} of $x$ in $X$ as follows
		\[
				w_X(x) = \frac{1 - q^{-(n+1)}}{1 - q^{-1}}  \frac{\max(1, \norm{x}^n)}{ \Nr(X,x)}.
		\]
		\end{definition}
		
		\begin{remark}\label{rem:wightFunction}
		    \begin{enumerate}
		        \item Notice that, on the $\Ocal$-points $X \cap \Ocal^N$ of the manifold $X$, we have $\Nr(X,x) = 1$ so the weight function $w_X$ is constant and takes the value
    			\[
    				w_X(x) = \frac{1 - q^{-(n+1)}}{1 - q^{-1}}, \quad x \in X \cap \Ocal^N.
    			\]
		        
		        \item The weight function $w_X$ is not intrinsic. As we shall see in \Cref{prop:wightfunction}, it depends on the probability distribution of the random element $(\bm{A},\bm{b})$ in $K^{n \times N} \times K^n$ is sampled~from.
		    \end{enumerate}
		\end{remark}
		           
		 \subsection{A pinch of intersection theory}\label{subsec:intersection}
		            
                 In this section we recall some facts from intersection theory of algebraic varieties (see \cite[Section 18]{harrisFirstCourse} and \cite{3264Eisenbud} for more details). The reader may skip this and come back to it when necessary.

				 Let $X \subset \A^N $ be an affine algebraic manifold of dimension $n$ and let $d$ be its degree. Then there exists a variety $\Vcal_{X}$ in $K^{n \times N } \times K^n$ of lower dimension such that
				 \begin{equation}\label{eq:intersectionAffine}
				    \#(\Lcal_{A,b} \cap X)  \leq d \quad \text{for any } (A,b) \in (K^{n \times N} \times K^n) \setminus \Vcal_{X}. 
				 \end{equation}
				 Since $\Vcal_X$ is a lower dimensional variety in $K^{n \times N} \times K^n$ it has measure $0$ with respect to the volume measure $\mathrm{d}A  \mathrm{d}b$. So if $(\bm{A}, \bm{b})$ is a random variable that has a density with respect to $\mathrm{d}A \mathrm{d}b$, then with probability $1$, the intersection $\Lcal_{\bm{A},\bm{b}} \cap X$ contains at most $d$ points. So, given a real valued function $f \colon X \to \R$ the function
				\[
				        \overline{f} \colon (A,b) \mapsto \sum_{x \in \Lcal_{A,b} \cap X} w_X(x) f(x),
				\]
				is well defined on the Zariski open set $K^{n \times N} \times K^n \setminus \Vcal_{X}$.

				Similarly, if $X \subset \P^{N-1}$ is a projective algebraic $K$-manifold of dimension $n$ and degree $d$ the set of matrices $A \in K^{n \times N}$ such that the intersection $\Lcal_{A} \cap X$ is infinite, where 
				\[
				    \Lcal_A = \{x \in \P^{N-1} \colon A x = 0\},
				\]
				is a lower dimensional algebraic variety $\mathcal{W}_X$ in $K^{n \times N}$. Hence $\mathcal{W}_X$ has measure $0$ with respect to the Haar measure $\mathrm{d}A$ and
				\begin{equation}\label{eq:intersectionProjective}
		            \# ( \Lcal_A \cap X ) \leq d \quad \text{for any } A \in K^{n \times N} \setminus \mathcal{W}_X.
				\end{equation}
				Moreover, given a real valued function $f \colon X \to \R$ the function
				\[
				        \overline{f} \colon A \mapsto \sum_{x \in \Lcal_{A} \cap X} f(x)
				\]
				is well defined on the Zariski open set $K^{n \times N} \setminus \mathcal{W}_X$.

		\subsection{The \texorpdfstring{$p$}~-adic co-area formula}
		In this section we recall a few notions on $p$-adic integration on manifolds. We refer the reader to \cite{BKL22, kulkarni2021integGeo,  mpopaInteg,LieGroups} and references therein  for a more detailed account.
		
		Let $X$ be a smooth algebraic (affine or projective) manifold defined over $K$. One can then endow the variety $X$ with the structure of a $K$-analytic manifold in the sense of \cite{BKL22} and a volume measure $\mu_X$. A definition of the latter is given in \Cref{eq:defOfAffineMeasure} for the affine case and in \Cref{eq:defOfProjMeasure} for the projective case.
		
		\begin{definition}
			Let $X$ and $Y$ be two $K$-analytic manifolds, $x \in X$ and $\varphi: X \to Y$ be a $K$-analytic map. We define the \emph{absolute Jacobian} of $\varphi$ at $x$ as 
			\[
			\J(\varphi,x) \coloneqq \N(\D_x \varphi),
			\]
			the absolute determinant of the differential $\D_x \varphi \colon T_x X \to T_{\varphi(x)} Y$ of $\varphi$ at the point $x$.
		\end{definition}
		
		The following is the $p$-adic coarea formula \cite[Theorem 3.3.2]{BKL22}.
		
		\begin{theorem}[Coarea formula] \label{thm:coarea}
			Let $X$ and $Y$ be two analytic $K$-manifolds such that $\dim(X) \geq \dim(Y)$ and let $\varphi: X \to Y$ be a $K$-analytic map. Then, for any function $f  \colon X \to \R$ that is integrable with respect to the volume measure on $X$, we have
			\[
			\int_{X}  \J(\varphi,x) f(x) \mu_X(\mathrm{d}x) = \int_{Y} \left( \int_{\varphi^{-1}(y)} f(z) \mu_{\varphi^{-1}(y)}(\mathrm{d}z) \right) \mu_Y(\mathrm{d}y).
			\]
		\end{theorem}

		\begin{corollary}\label{cor:densities}
			Let $X$ and $Y$ be two $K$-manifolds and $\varphi \colon X \to Y$ an analytic map from $X$ to $Y$.
			
			\begin{enumerate}[label=(\roman*)]
				\item \label{cor:densitiespart1} Suppose that $\bm{\xi}$ is an $X$-valued random variable with density $f$ with respect to $\mu_X$. Then the density $g$ of $\bm{\eta} = \varphi(\bm{\xi})$ with respect to $\mu_Y$ is
				\[
				g(y) = \int_{\varphi^{-1}(y)}   \frac{f(z)}{\J(\varphi, z)}  \mu_{\varphi^{-1}(y)}(\mathrm{d}z).
				\]
				
				\item \label{cor:densitiespart2} Let $\bm{\eta}$ be a $Y$-valued random variable with density $g$ with respect to $\mu_Y$ and let $\bm{\xi}$ be the $X$-valued random variable such that, conditioned on $(Y=y)$, the variable $\bm{\xi}$ has density $f_{y}$ on $\varphi^{-1}(Y)$ with respect to $\mu_{\varphi^{-1}(y)}$. Then $\bm{\xi}$ has density
				\[
				f(x) = \J(\varphi, x) g(\varphi(x)) f_{\varphi(x)}(x).
				\] 
			\end{enumerate}
		\end{corollary}
		
		\begin{proof}
			\begin{enumerate}[label=(\roman*)]
				\item Let $V$ be a Borel set in $Y$. Applying \Cref{thm:coarea} we get
				\begin{align*}
					P(\bm{\eta} \in V) &= \int_X 1_{V}(\varphi(x)) f(x) \mu_X(\mathrm{d}x) \\
					&= \int_Y \left( \int_{\varphi^{-1}(y)} \frac{f(z)}{\J(\varphi, z)} 1_{V}(\varphi(z))  \mu_{\varphi^{-1}(y)} (\mathrm{d}z) \right)  \mu_Y(\mathrm{d}y) \\
					&= \int_Y \left( \int_{\varphi^{-1}(y)} \frac{f(z)}{\J(\varphi, z)}  \mu_{\varphi^{-1}(y)} (\mathrm{d}z) \right) 1_{V}(y)  \mu_Y(\mathrm{d}y) \\
					&= \int_Y g(y) \mu_Y(\mathrm{d}y).
				\end{align*}
					
				\item  	Let $U$ be a Borel set in $X$. Then, applying \Cref{thm:coarea} we get
				\begin{align*}
					P(\bm{\xi} \in U) &= \E[ P(\bm{\xi} \in U  \mid \bm{\eta})]\\
					&= \int_{Y}  \left( \int_{\varphi^{-1}(y)}  f_{y}(z)  1_U(z)  \mu_{\varphi^{-1}(y)}(\mathrm{d}z) \right) g(y) \mu_Y(\mathrm{d}y) \\
					&= \int_{Y}  \left( \int_{\varphi^{-1}(y)}  g(\varphi(z)) f_{\varphi(z)}(z)  1_U(z)  \mu_{\varphi^{-1}(y)}(\mathrm{d}z) \right) \mu_Y(\mathrm{d}y) \\
					&= \int_{X}  \J(\varphi, x) g(\varphi(x)) f_{\varphi(x)}(x)  1_U(x)   \mu_X(\mathrm{d}x) \\
					&= \int_{X}  f(x) 1_U(x)  \mu_X(\mathrm{d}x).  
				\end{align*}
			\end{enumerate}
		\end{proof}

		We denote by $\Gr(n,K^m)$ the Grassmannian variety parametrizing $n$-dimensional vector subspaces of $K^m$. The orthogonal group $\GL(m,\Ocal)$ has a natural action on $\Gr(n,K^m)$.
		
		\begin{lemma}\label{lem:uniqueOrthogonalProb}
		    Let $m \geq n \geq 1$ be two integers. There exists is a unique orthogonally invariant probability distribution on the Grassmanian $\Gr(n, K^{m})$.
		\end{lemma}
		\begin{proof}
		    Since $\GL(m,\Ocal)$ acts transitively on $\Gr(n,K^m)$ and the stabilizer of the subspace generated by the first $n$ vectors of the standard basis of $K^m$ is
		    \[
		         H = \left \{ \begin{pmatrix}
                A & \rvline & C  \\ 
                   \hline
                0 & \rvline & B \end{pmatrix} \colon A \in \GL(n, \Ocal), B \in \GL(m-n,\Ocal) \quad \text{and} \quad C \in \Ocal^{n \times (m-n)} \right  \},
		    \]
		    we can write $\Gr(n,K^m)$ as a homogeneous space as follows
		    \[
		        \Gr(n,K^m) = \GL(m,\Ocal) / H.
		    \]
		    Let $\nu$ be a probability measure on $\GL(m,\Ocal) / H$ that is $\GL(m,\Ocal)$-invariant. Then its pullback $\nu^\ast$ to $\GL(m,\Ocal)$ is also $\GL(m,\Ocal)$ invariant, so it is a Haar measure on $\GL(m,\Ocal)$ with $\nu^{\ast}(\GL(m,\Ocal)) = 1$. We then conclude since $\GL(m,\Ocal)$ is a compact topological group, there is a unique Haar measure on $\GL(m, \Ocal)$ up to scaling.
		\end{proof}
		
        We end this section with the following simple but useful lemma.
        
        \begin{lemma}\label{lem:ExpectationDet}
		Let $n \geq 1 $ be a positive integer, then we have
		\[
				\int_{C \in \Ocal^{n \times n}} |\mathrm{det}(C)| \ \mathrm{d}C = \frac{1 - q^{-1}}{1 - q^{-(n+1)}}.
		\]
		\end{lemma}
		\begin{proof}
			Let $\bm{C}$ be a random matrix in $K^{n \times n}$ whose entries are independent and uniform in $\Ocal$. Then the integral in question is the expectation $\E[ |\mathrm{det}(\bm{C})|]$. We can compute this expectation using the distribution of $|\mathrm{det}(\bm{C})|$ from \cite[Theorem 4.1]{ElementaryDivisors}. We then have
			\[
				\E[|\mathrm{det}(\bm{C})|] =  (1 - q^{-1}) \cdots (1 - q^{-n}) \sum_{m  = 0}^{\infty} \binom{n +m -1}{m}_{q^{-1}} q^{- 2m},
			\]
            where $\binom{n}{k}_{q^{-1}}$ denotes the usual $q^{-1}$-binomial coefficient (also known as the Gaussian binomial coefficient):
                \[
                \binom{n}{k}_{q^{-1}} \coloneqq \frac{(1 - q^{-1}) \cdots (1 - q^{-n})}{ (1 - q^{-1}) \cdots (1 - q^{-k}) \times (1 - q^{-1}) \cdots (1 - q^{-(n-k)}) } \quad \text{for } n \geq k \geq 0.
                \]
			Then using the well known generating series
			\[
				\sum_{m  = 0}^{\infty} \binom{n +m -1}{m}_{q^{-1}} t^{m} = \prod_{k = 0}^{n-1} \frac{1}{1 - q^{-k}t},
			\]
			we get
            $$ \sum_{m  = 0}^{\infty} \binom{n +m -1}{m}_{q^{-1}} q^{-2m} = \prod_{k = 0}^{n-1} \frac{1}{1 - q^{- k - 2 }}  = \frac{1}{(1-q^{-2}) \cdots (1- q^{-(n+1)})}.$$
            So we deduce that 
            \[
                \E[|\det(\bm{C})|] = \int_{C \in \Ocal^{n \times n}} |\det(C)| \ \mathrm{d}C = \frac{1 - q^{-1}}{1 - q^{-(n+1)}}.  
            \]
		\end{proof}
	
		\begin{remark}
		    We shall see in \Cref{sec:Sampling_ProjectiveManifolds_Prelim} that the integral in \Cref{lem:ExpectationDet} is exactly the total measure projective pace $\mu_{\P^n}(\P^n)$ i.e.
		    $$\mu_{\P^n}(\P^n) = \frac{1 - q^{- (n+1)}}{1 - q^{-1}}.$$
		\end{remark}

    \section{Sampling from affine manifolds}\label{sec:Sampling_AffineManifolds}

		In this section, we proceed to proving the main results of this article, namely Theorems \ref{thm:integral_EP} and \ref{thm:xi_density}. Similar results for projective manifolds are stated and proved in \Cref{sec:Sampling_ProjectiveManifolds}. We start with the following:
        
		\begin{proposition}\label{prop:wightfunction} 
    		Let $X \subset \A^N$ be an affine algebraic $K$-manifold of dimension $n$ and $x$ a point on $X$. We then have
    		\[
    		w_X(x) = \left( \int_{A \in K^{n \times N} }  | \mathrm{det}(A_{|T_{x}X}) | \ 1_{ A \in \Ocal^{n \times N} ,  \norm{Ax} \leq 1} \mathrm{d}A \right) ^{-1}\ .
    		\]
		\end{proposition}
		\begin{proof}
        		Let $U \in \GL(N,\Ocal)$ such that $y \coloneqq  U x = (0, \dots, 0,\varpi^{\val(x)})^{\top}$. Let $W$ be a matrix whose columns form an orthonormal basis of $T_{x}X$. Let $R_x$ denote the $N\times N$ matrix $R_x = \diag(1, \dots, 1, \varpi^{v(x)})$. We denote by $I_X(x)$ the following integral:
        		\[
        		I_X(x) = \int_{A \in K^{n \times N}} |\mathrm{det}(A_{|T_{x}X})| 1_{A \in \Ocal^{n \times N}} 1_{ Ax \in \Ocal^n} \mathrm{d}A.
        		\]
        		 Then, by a change of variables $BU = A$, we have
        		\begin{align*}
        		I_X(x) &= \int_{B \in K^{n \times N}} |\mathrm{det}((BU)_{|T_{x}X})| \ 1_{B \in \Ocal^{n \times N}} 1_{ By \in \Ocal^n } \mathrm{d}B\\
  				       &= \int_{B \in K^{n \times N}} |\mathrm{det}(BUW)|  \  1_{B \in \Ocal^{n \times N}} 1_{ B R_x \in \Ocal^{n \times N}} \mathrm{d}B\\
				       &=\int_{B \in  \Ocal^{n \times N} \cap ( \Ocal^{n \times N} R_x^{-1} ) }  |\mathrm{det}(BUW)| \mathrm{d}B.
        		\end{align*}
        		Notice the following equality:
        		\[
        		\Ocal^{n \times N} \cap \Ocal^{n \times N}R_x^{-1} = \Ocal^{n \times N} S_x,
        		\]
        		where $S_x = \diag(1,\dots, 1, \varpi^{\max(0,-\val(x))}) \in K^{N \times N}$. So, using the change of variables $B = B' S_x$, we deduce that
        		\begin{align*}
        						 I_X(x) 	&= \int_{B \in \Ocal^{n \times N}S_x} |\mathrm{det}(BUW)| \mathrm{d}B\\
        		         						&= \left(\frac{1}{\max(1, \norm{x})}\right)^{n} \int_{B' \in \Ocal^{n \times N}} | \mathrm{det}(B' S_x U W)| \mathrm{d}B'.
        		\end{align*}
        		Let us write the Smith normal form of the matrix $S_x U W$, i.e. 
        		\[
        				S_x U W = V_1 D V_2,
        		\] 
        		where $V_1 \in \GL(N,\Ocal), V_2 \in \GL(n,\Ocal)$ and $D = \diag(\varpi^{v_1}, \dots, \varpi^{v_n}) \in K^{N \times n}$. So, by the change of variables $B'V_1 = C$, we get
        		\begin{align*}
        								I_X(x) &= \frac{1}{\max(1, \norm{x}^n)}  \int_{B' \in \Ocal^{n \times N}} | \mathrm{det}(B' V_1 D V_2)| \mathrm{d}B'\\
        										   &= \frac{1}{\max(1,\norm{x}^n)} \int_{C \in \Ocal^{n \times N}}  |\mathrm{det}(C D)| \mathrm{d}C \\
        										   &= \frac{q^{-(v_1 + \cdots + v_n)}}{\max(1,\norm{x}^n)} \int_{C \in \Ocal^{n \times n}}  |\mathrm{det}(C)| \mathrm{d}C.
        		\end{align*}
        		Combining the previous equation with \Cref{def:NrDefAffine}, \Cref{def:affineWeight} and \Cref{lem:ExpectationDet}, we get
        		\[
        			I_X(x) = \frac{\Nr(X,x)}{\max(1, \norm{x}^n)}  \frac{1 - q^{-1}}{1 - q^{-(n+1)}}  = \frac{1}{w_X(x)}
        		\] 
        		as desired.
		\end{proof}
  
		\begin{remark}\label{rem:weightBehaviour}
		    \begin{enumerate}[label=(\roman*)]
		        \item \Cref{prop:wightfunction} gives another proof of the fact that $\Nr(X,x)$ does not depend on the choice of $U$ and $W$ in \Cref{def:NrDefAffine}.
		        \item Recall, from \Cref{rem:wightFunction}, that the weight function is constant on $X\cap \Ocal^N$. Unwinding the definition of $w_X$ we can also see that for $U \in \GL(N,\Ocal)$ and $x \in X$ we have $w_{U X}(U x) = w_{X}(x)$ where $UX = \{U x \colon x \in X\}$.
		        \item \label{rem:tip}  If the probability density $f$ we wish to sample from is supported on $X \cap \varpi^{-r} \Ocal^N$, we can scale $X$ by $\varpi^r$ and sample $\bm{\xi'}$ from $\varpi^{r} X \cap \Ocal^N$ (where the weight function is constant) with density $f(\varpi^{-r} \cdot)$. We can then obtain a random variable $\bm{\xi}$ on $X$ with density $f$ by taking $\bm{\xi} = \varpi^{-r} \bm{\xi}'$.
		    \end{enumerate}
		\end{remark}
		
		We are now ready to prove our main theorems.
		
		\begin{proof}[Proof of \Cref{thm:integral_EP}]

		By definition we have
		\[
			\E( \overline{f}(\bm{A},\bm{b})) = \int_{K^{n \times N} \times K^{n}} \overline{f}(A,b) 1_{A \in \Ocal^{n \times N}, b \in \Ocal^n} \mathrm{d}A \mathrm{d}b.
		\]
		Let us define the following map:
		\begin{equation} \label{eq:varphiAffine}
			 \varphi \colon K^{n \times N} \times X \to K^{n \times N} \times K^{n} , \quad  (A,x) \mapsto (A,Ax).
		\end{equation}
		The map $\varphi$ is analytic and its differential is given by 
		\[
			D_{(A,x)} \varphi \colon K^{n \times N} \times T_{x}X \to K^{n \times N} \times K^n, \quad   (H,u) \mapsto (H, Hx + Au).
		\]
		So the differential of $\varphi$ at $(A,x)$ acts trivially on the first component of the product $K^{n \times N} \times T_xX$ and acts as $A_{|T_xX}$ on the second, i.e. its determinant is given by
		\begin{equation} \label{eq:JacAndDet}
			\J(\varphi, (A,x)) = \left| \mathrm{det}(A_{|T_xX})  \right|.
		\end{equation}
		 Applying \Cref{thm:coarea} for the function $\varphi$ yields
		 	   \begin{align*}
						 & \int_{K^{n \times N} \times X}  \left| \mathrm{det}(A_{|T_xX}) \right|   w_X(x) f(x) 1_{A \in \Ocal^{n \times N}, \ Ax \in \Ocal^n}  \ \mathrm{d}A \ \mu_X(\mathrm{d}x)\\
						 = &  \int_{K^{n \times N} \times X}  \J(\varphi, (A,x)) w_X(x) f(x)  1_{A \in \Ocal^{n \times N}, \ Ax \in \Ocal^n} \ \mathrm{d}A \ \mu_X(\mathrm{d}x) \\
						 = & \int_{K^{n \times N} \times K^{n}} \left( \int_{\varphi^{-1}(A,y)}  w_X(z)f(z) 1_{A \in \Ocal^{n \times N}, \ Az \in \Ocal^n } \ \mu_{\varphi^{-1}(A,y)}(\mathrm{d}z) \right) \mathrm{d}A \mathrm{d}y\\
						 = & \int_{K^{n \times N} \times K^{n}} \left( \int_{\varphi^{-1}(A,y)}  w_X(z)f(z)   \mu_{\varphi^{-1}(A,y)}(\mathrm{d}z) \right) 1_{A \in \Ocal^{n \times N}, \ y \in \Ocal^n } \ \mathrm{d}A \mathrm{d}y.
			  \end{align*}
		 But, for $A \in K^{n \times N}$ and $y \in K^n$ we have
		 \[
		 	\varphi^{-1}((A,y)) = \{ (A,z) \in K^{n \times N} \times X \colon Az = y  \text{ and } z \in X \},
		 \]
		 and this is a finite set for almost every $A$ and $y$. So for almost every $A$ and $y$, the measure $\mu_{\varphi^{-1}(A,y)}$ equals the counting measure on the finite set $\varphi^{-1}((A,y))$ and we then have
				 \begin{align*}
				 &\int_{K^{n \times N} \times X}  |\mathrm{det}(A_{|T_xX})|   w_X(x) f(x) 1_{A \in \Ocal^{n \times N}, \ Ax \in \Ocal^n}  \ \mathrm{d}A  \mu_X(\mathrm{d}x) \\
		   = & \int_{K^{n \times N} \times K^{n}}  \left( \sum_{\substack{x \in X, \\  Ax = y} } w_X(x) f(x)  \right)  1_{A \in \Ocal^{n \times N}, y \in \Ocal^n } \ \mathrm{d}A \mathrm{d}y \\
		   = & \int_{K^{n \times N} \times K^{n}} \overline{f}(A,y) 1_{A \in \Ocal^{n \times N}, \ y \in \Ocal^n } \ \mathrm{d}A \mathrm{d}y\\
		    =& \  \E[\overline{f}(\bm{A},\bm{b})].
				\end{align*}
		Hence the equation
		\begin{equation}\label{eq:lastProof1}
			\E \left[\overline{f}(\bm{A},\bm{b}) \right] =   \int_{X} \left( \int_{K^{n \times N}} |\mathrm{det}(A_{|T_xX})| 1_{A \in \Ocal^{n \times N}, \ Ax \in \Ocal^n} \mathrm{d}A \right)  w_X(x) f(x) \mu_X(\mathrm{d}x).
		\end{equation}
		Then, combining \Cref{eq:lastProof1} and \Cref{prop:wightfunction}, we conclude that
		\[
			\E\left[\overline{f}(\bm{A},\bm{b})\right] =   \int_{X} f(x) \mu_X(\mathrm{d}x).  
		\]
		\end{proof}
		\begin{proof}[Proof of \Cref{thm:xi_density}]
    		Let $(\bm{\widetilde{A}}, \bm{\xi})$ be the random variable, with values in $K^{n \times N} \times X$, obtained by first sampling $(\bm{\widetilde{A}}, \bm{\widetilde{b}}) \in K^{n \times N} \times K^n$ with distribution $\overline{f}(A,b) 1_{A \in \Ocal^{n \times N}, b \in \Ocal^n} \mathrm{d}A \mathrm{d}b$ and then choosing a point $\bm{\xi}$ from $\Lcal_{ \bm{\widetilde{A}}, \bm{\widetilde{b}} } \cap X$ with probability 
    		\[\frac{w_X(x) f(x)}{\overline{f}(\bm{\widetilde{A}}, \bm{\widetilde{b}})}.\] 
    		Then applying \Cref{cor:densities}-\ref{cor:densitiespart2} to the map $\varphi$ from \Cref{eq:varphiAffine}, we deduce that $(\bm{\widetilde{A}}, \bm{\xi})$ has density
    		\begin{align*}
    				g_{(\bm{\widetilde{A}}, \bm{\xi})} (A,x) &= \overline{f}\left( \varphi(A,x) \right) 1_{\varphi(A,x) \in \Ocal^{n \times N} \times \Ocal^n} \frac{w_X(x)f(x)}{\overline{f}(\varphi(A,x))} \J(\varphi, (A,x))\\
    																		 &= w_X(x) f(x)   1_{\varphi(A,x) \in \Ocal^{n \times N} \times \Ocal^n} \J(\varphi, (A,x))
    		\end{align*}
    		with respect to the volume measure $\mathrm{d}A \mu_X(\mathrm{d}x)$ on $K^{n \times N} \times X$. Computing the second marginal of this joint distribution, we deduce that the density $g_{\bm{\xi}}$ of $\bm{\xi}$ is
    		\begin{align*}
    		g_{\bm{\xi}}(x) &= \int_{A \in K^{n \times N} } w_X(x) f(x) \J(\varphi, (A,x)) 1_{\varphi(A,x) \in \Ocal^{n \times N} \times \Ocal^n} \mathrm{d}A\\
    					  &=  w_X(x) f(x)  \int_{A \in K^{n \times N} } |\mathrm{det}(A_{|T_x X})| 1_{A \in \Ocal^{n \times N}, Ax \in \Ocal^n} \mathrm{d}A\\
    					  & = f(x).
    		\end{align*}
    		The second (resp. third) equation follows from \Cref{eq:JacAndDet} (resp. \Cref{prop:wightfunction}). So, as desired, $\bm{\xi}$ has density $f$ with respect to $\mu_X$ on $X$. 
		\end{proof}

\section{Sampling from projective manifolds}\label{sec:Sampling_ProjectiveManifolds}

		This section deals with sampling from projective algebraic manifolds. More precisely, we shall state and prove analogs of \Cref{thm:integral_EP} and \Cref{thm:xi_density} in projective space.

		Let $N \geq 2$ be an integer. We denote by $\P^{N-1}$ the projective space of dimension $N-1$ over $K$. Let us denote by $S^{N-1}$ the unit sphere in $K^N$, i.e.,
		\[
		    S^{N-1} \coloneqq \{ x \in K^N \colon \norm{x} = 1 \}.
		\]
		We warn the reader that, unlike the Euclidean setting, the unit sphere is actually an open set in $K^N$ and has dimension $N$ (as a topological space). Consider the \emph{Hopf fibration}
		\[
		\psi \colon S^{N-1} \to \P^{N-1}, \quad (x_1, \dots, x_N) \mapsto (x_1 \colon \cdots \colon x_N).
		\]
		The projective space $\P^{N-1}$ can be endowed with a metric $d$ defined as follows:
		\[
		    d(x,y) = \norm{\tilde{x} \wedge \tilde{y}}, \quad x,y \in \P^{N-1}
		\]
		where $\tilde{x} \in \psi^{-1}(x)$, $\tilde{y} \in \psi^{-1}(y)$ and the norm $\norm{\tilde{x} \wedge \tilde{y}}$ is the standard norm in $\bigwedge\limits^2 K^N$ associated to its standard lattice $\bigwedge\limits^2 \Ocal^N$. This metric is called the \textit{Fubini-Study} metric. 	For $x \in \P^{N-1}$ and $\epsilon > 0$ let us denote by 
		\[\B_{N-1}(x,\epsilon) \coloneqq \{ y \in \P^{N-1} \colon d(x,y) \leq \epsilon \}\]
		the ball of radius $\epsilon$ around $x$.

		Endowed with the metric $d$, the projective space $\P^{N-1}$ is a compact metric space on which we define a volume measure $\mu_{\P^{N-1}}$ as follows
		\[
		    \mu_{\P^{N-1}} \coloneqq \frac{1}{1 - q^{-1}} \psi_\ast \mu_{S^{N-1}},
		\]
		that is the normalized push-forward of $\mu_{S^{N-1}}$ by $\psi$. Note that $\mu_{S^{N-1}}(S^{N-1}) = 1 - q^{-N}$, so the measure $\mu_{\P^{N-1}}$ is finite and we have
		\[
		    \mu_{\P^{N-1}}(\P^{N-1}) = \frac{1 - q^{N}}{1 - q^{-1}}.
		\]

        \begin{remark}
            Notice that from \Cref{lem:ExpectationDet}, we have
            \[
                \int_{C \in \Ocal^{n \times n}} |\det(C)| \mathrm{d}C = \frac{1}{\mu_{\P^n}(\P^n)}.
            \]
        \end{remark}

		A projective algebraic variety in $\P^{N-1}$ is the zero set of a system of homogeneous polynomials $\bm{p} = (p_1, \dots, p_r)$ in $K[x_1, \dots, x_N]$; that is
		\[
		    \{  x \in \P^{N-1} \colon p_1(x) = \cdots = p_r(x) = 0\}.
		\]
		We refer to irreducible and smooth projective varieties as projective algebraic manifolds. 
		
		Let $X \subset \P^{N-1}$ be an algebraic projective manifold of dimension $n \geq 1$. Similar to the affine case (\ref{eq:defOfAffineMeasure}), we can define a volume measure on $X$ as follows:
		\begin{equation}\label{eq:defOfProjMeasure}
			\mu_X(V) \coloneqq  \lim\limits_{\epsilon \to 0} \frac{  \mu_{\P^{N-1}} \left(  \bigcup\limits_{ x \in V} \B_{N-1}(x, \epsilon)  \right) }{ \mu_{\P^{N-1-n}} \left(  \B_{N-1-n}(0,\epsilon) \right) } , \quad \quad  \text{for } V \subset X \text{ open}.
		\end{equation}
		The limit in (\ref{eq:defOfProjMeasure}) exists (see \cite{Serre81,kulkarni2021integGeo} for more details) and this defines a volume measure $\mu_X$ on the projective manifold $X$.
		
		\begin{remark}
		    For our purposes, the main difference between the affine and projective spaces is that the projective space is a compact topological space (with the quotient topology induced by the Hopf fibration $\psi$). So, unlike the affine case, a projective algebraic manifold admits a uniform probability density. Also, loosely speaking, there are no ``far'' points in the projective space, so as we shall see, the weight function is constant or, in other words, no point gets more weight than another. We can say that the space is, in some sense, ``isotropic".
		\end{remark}
		
		Before we state our results for projective manifolds, we recall a few facts and establish a couple of preliminary results.
		
		\subsection{Preliminaries} \label{sec:Sampling_ProjectiveManifolds_Prelim}
		
		Suppose that $X \subset \P^{N-1}$ is a projective algebraic manifold of dimension $n$ defined by homogeneous polynomials $p_1,\dots, p_r \in K[x_1, \dots, x_N]$ and let $x$ be a point in $X$.
		The tangent space $T_x X$ can be defined in many ways, and one way to do so is the following. The cone $\widetilde{X} \subset \A^{N}$ over $X$ defined as follows
		\[
		    \widetilde{X} = \{  (\lambda y_1, \dots, \lambda y_N)\in \A^{N} \colon \lambda \in K \text{ and } (y_1 : \dots : y_N) \in X  \}.
		\]
		This is an affine algebraic variety which is smooth at every non-zero point $x \in \widetilde{X} \setminus \{0\}$ and has dimension $n+1$. The tangent space $T_x \widetilde{X}$ is a linear subspace in $K^N$ of dimension $n+1$ and $x \in T_x \widetilde{X}$. The tangent space $T_xX$ can then be defined as an orthogonal complement\footnote{All such vector spaces are isomorphic to one another.} of the line $K \cdot x$ in $T_x \widetilde{X}$ and we thus view $T_xX$ as a linear subspace\footnote{The projective tangent space is also often defined as the projectivization of $T_x \widetilde{X}$.} of $K^N$ of dimension $n$.

		\begin{proposition}\label{prop:TangentSpace}
			Let $X \subset \P^{N-1}$ be a projective algebraic manifold of dimension $n$ and let us define $\Xcal \subset \A^{n \times N} \times \P^{N-1}$ as follows:
			\[
				\Xcal = \{ (A,x) \in  \A^{n \times N} \times X \colon A x = 0 \}.
			\]
			Then $\Xcal$ is a manifold, and for $(A,x) \in \Xcal$ we have
			\[
				T_{(A,x)} \Xcal = \{ (H,h) \in K^{n \times N} \times T_{x}X \colon Hx + Ah = 0 \}.
			\]
			Moreover, if $\varphi, \phi$ are the projections from $\Xcal$ to $\A^{n \times N}$ and $\P^{N-1}$ respectively, then we have
			\[
				\frac{\J(\varphi,(A,x))}{\J(\phi, (A,x))} = |\mathrm{det}(A_{|T_x X})|,
			\]
			for $(A,x) \in \Xcal$ such that $A_{|T_xX}$ is an isomorphism.
		\end{proposition}
		\begin{proof}
		Let $(p_1, \dots, p_r) \in K[x_1, \dots, x_{N}]$ be homogeneous polynomials generating the ideal of $X$. Let $(A,x) \in \Xcal$ and let $J_x$ be the following Jacobian matrix
		\[
			J_x = \left( \frac{\partial p_i}{\partial x_j} (x) \right)_{1 \leq i \leq r, 1 \leq j \leq N}.
		\]
		Then, considering $\Xcal$ as the variety in $\A^{n \times N} \times \P^{N-1}$ cut out by the equations $A x = 0 $ and $p_1(x) = \cdots = p_r(x) = 0$ we can compute the Jacobian matrix of $\Xcal$ at the point $(A,x)$. This matrix represents the linear map
		\[
		K^{n \times N} \times K^N \to  K^{n} \times K^{r}, \quad (H, h) \mapsto (Hx + Ah, J_x h).
		\]
		The tangent space of $\Xcal$ at $(A,x)$ is the kernel of this map, so
		\[
		T_{(A,x)} \Xcal = \{ (H,h) \in K^{n \times N} \times T_{x}X \colon Hx + Ah = 0 \}.
		\]
		The projection maps $\varphi, \phi$ are clearly analytic, and for any $(A,x) \in \Xcal$ we have
        \begin{align*}
        d_{(A,x)} \varphi \colon T_{(A,x)} \Xcal&\to K^{n \times N} & d_{(A,x)} \phi \colon T_{(A,x)} \Xcal&\to T_xX\\
        (H,h)&\mapsto H& (H,h)&\mapsto h.
        \end{align*}
        
		Suppose that $(A,x) \in \Xcal$ is such that $A_{|T_xX}$ is an isomorphism. Fix $U \in \GL(N,\Ocal)$ such that $Ux = (1:0:\cdots:0)^\top$  and define the maps
        \begin{align*}
        \pi_1  \colon K^{n \times N} &\to T_{(A,x)} \Xcal & \pi_2  \colon T_x X & \to T_{(A,x)} \Xcal \\
        H & \mapsto (H, - (A_{| T_xX})^{-1} H x ) & h & \mapsto ((-Ah | 0)U,  h )
        \end{align*}
		where $(-Ah|0) \in K^{n \times N}$. Notice that $d_{(A,x)}\varphi \circ \pi_1 = \mathrm{Id}_{K^{n \times N}}$ and $d_{(A,x)}\phi \circ  \pi_2  = \mathrm{Id}_{T_x X}$. Since $A \in \Ocal^{n \times N}$ we have
		\[
				 \N(\pi_2) = 1,
		\]
		because $\pi_2$ sends any orthonormal basis of $T_x X$ to an orthonormal family in $T_{(A,x)} \Xcal \subset K^{n \times N} \times K^N$. Also, since $A \in \Ocal^{N \times N}$, the singular values of $A_{|T_{x}X}$ are all in $\Ocal$ so the singular values of $A_{|T_x X}^{-1}$ have negative or zero valuation. From this we can see that
		\[
				\N(\pi_1) = |\mathrm{det}(A_{|T_x X})|^{-1}.
		\]
		We deduce that
        \[
                \frac{\J(\varphi, (A,x))}{ \J(\phi,(A,x))} =  \frac{\N(\pi_2)}{ \N(\pi_1) } = |\mathrm{det}(A_{|T_x X})|.  
        \]
		\end{proof}
		
		\begin{lemma}\label{lem:projectiveWeight}
			Let $X$ be a projective manifold of dimension $n$ in $\P^{N-1}$ and $x$ be a point on $X$. Set $M_x \coloneqq \{ A \in K^{n \times N} \colon A x = 0 \}$. Then
			\[
			\int_{M_x} | \mathrm{det}(A_{|T_xX} )|   1_{A \in \Ocal^{n \times N}} \mu_{M_x}(\mathrm{d}A) = \frac{1 - q^{-1}}{1 - q^{-(n+1)}}.
			\]
		\end{lemma}
		\begin{proof}
		    Let $W \in \Ocal^{N \times n}$ be a matrix whose columns form an orthonormal basis of $T_x X$ and let $U \in \GL(N,\Ocal)$ such that $U x = e_1 = (1: 0 : \cdots : 0)^\top$. The space $M_x$ is a vector space of dimension $(N-1) \times n$ and $M_{e_1} U = M_x$. So, with the change of variable $BU = A$, we get
			\begin{align*}
        			\int_{M_x} | \mathrm{det}(A_{|T_xX} )| 1_{A \in \Ocal^{n \times N}} \mu_{M_x}(\mathrm{d}A) 
                    &= \int_{A \in M_x} |\mathrm{det}(AW)| 1_{A \in \Ocal^{n \times N}} \mu_{M_x}(\mathrm{d}A) \\
                    &= \int_{B \in M_{e_{1}}} |\mathrm{det}(BUW)| 1_{B \in \Ocal^{N \times n}}\ \mu_{M_x}(\mathrm{d}A)\\
                    &= \int_{C \in \Ocal^{n \times (N-1)}}  |\mathrm{det}( (0 | C) U W)| \ \mathrm{d}C.
		 	\end{align*}
			Let $\widetilde{W} \in K^{(N-1) \times n}$	be the matrix obtained from $U W$ by deleting the first row and let us write the Smith normal form of $\widetilde{W}$ as
			\[
				\widetilde{W} = V_1 D V_2,
			\]
			where $V_1 \in \GL(N-1, \Ocal), V_2 \in \GL(n,\Ocal)$ and $D = \diag(\varpi^{v_1}, \dots, \varpi^{v_n}) \in K^{(N-1) \times n}$. We then deduce that
			\begin{align*}
			\int_{M_x} | \mathrm{det}(A_{|T_xX} )| 1_{A \in \Ocal^{n \times N}} \mu_{M_x}(\mathrm{d}A) &= \int_{C \in \Ocal^{n \times (N-1)}}  |\mathrm{det}( C \widetilde{W})| \ \mathrm{d}C\\
									&=  q^{- (v_1 + \cdots + v_n) }  \int_{C \in \Ocal^{ n \times N}}  |\mathrm{det}(C)| \ \mathrm{d}C\\
								    &=  \N(\widetilde{W}) \frac{1 - q^{-1}}{1 - q^{-(n+1)}}.
			\end{align*}
            Since $X$ is a projective manifold, the tangent space $T_xX$ is orthogonal to $x$ (see \Cref{sec:Sampling_ProjectiveManifolds_Prelim}). We deduce that the columns of $UW$ are orthogonal to $(1,0, \dots, 0)^\top$ so $\widetilde{W}$ has orthonormal columns. Hence $\N(\widetilde{W}) = 1$ which finishes the proof.
	\end{proof}

		Similarly to the affine case given a real valued function $f \colon X \to \R$, we define the weighted average function of $f$ as follows:
		\[
			\overline{f}(A) = \sum_{x \in \Lcal_{A} \cap X} f(x), \quad \text{for } A \in K^{n \times N}.
		\]
		By convention, the sum is taken to be $0$ whenever $\Lcal_A \cap X$ is empty or infinite.

		\subsection{Sampling from projective manifolds}
		
		Now we state and prove the analogues of Theorems \ref{thm:integral_EP} and \ref{thm:xi_density} for the projective case.
		\begin{theorem}\label{thm:integral_EPprojective}
			Let $X \subset \P^{N-1}$ be an $n$-dimensional projective algebraic manifold defined over $K$. Let $\bm{A}$ be a random variable in $K^{n \times N}$ with distribution $1_{A \in \Ocal^{n \times N}} \mathrm{d}A$. Then we have
			\[
				\int_{X}  f(x) \mu_X(\mathrm{d}x) = \mu_{\P^n}(\P^n) \ \E[\overline{f}(\bm{A})].
			\]
		\end{theorem}
		\begin{proof}
			Let $\Xcal \subset \A^{n \times N} \times \P^{N-1}$ be the algebraic  variety defined by
			\[
				\Xcal \coloneqq  \left\{ (A,x) \in \A^{n \times N} \times X \colon  Ax = 0 \right\}.
			\]
			Let us denote by $\varphi$ and $\phi$ the natural projections from $\Xcal$ onto $K^{n \times N }$ and $X$ respectively, and, for a point $x\in X$, set $M_x \coloneqq \{ A \in K^{n \times N} \colon A x = 0 \}$. We apply \Cref{thm:coarea} and \Cref{prop:TangentSpace} on $\varphi$ and then on $\phi$ to get the following:
			\begin{align*}
			\E \left[ \overline{f} (\bm{A}) \right] &=  \int_{K^{n \times N}} \left( \sum_{\substack{x \in X, \\ Ax = 0}}   f(x) \right) 1_{A \in \Ocal^{n \times N}} \mathrm{d}A\\
												    &= \int_{K^{n \times N}} \left(\  \int_{(A,z) \in \varphi^{-1}(A)}   f(z)   1_{A \in \Ocal^{n \times N}}  \mu_{\varphi^{-1}(A)}(\mathrm{d}z) \right) \mathrm{d}A\\
												    &= \int_{\Xcal}  \J(\varphi,(A,x)) f(x)  1_{A \in \Ocal^{n \times N}} \mu_{\Xcal}(\mathrm{d}A,\mathrm{d}x)\\
												    &= \int_{X} \left( \int_{(A,x) \in \phi^{-1}(x)}   \frac{\J(\varphi,(A,x))}{\J(\phi,(A,x))} 1_{A \in \Ocal^{n \times N}} \mu_{\phi^{-1}(x)} \mathrm{d}A \right)f(x) \mu_X(\mathrm{d}x)\\
  												    &= \int_{X}  \left( \int_{A \in M_x} |\mathrm{det}(A_{|T_x X})| 1_{A \in \Ocal^{n \times N}}  \mu_{\phi^{-1}(x)}(\mathrm{d}A) \right) f(x) \mu_X(\mathrm{d}x)\\
  												    &= \frac{1 - q^{-1}}{1 - q^{-(n+1)}} \int_X f(x) \mu_X(\mathrm{d}x)\\
                                                    &= \frac{1}{\mu_{\P^n}(\P^n)} \int_X f(x) \mu_X(\mathrm{d}x).
			\end{align*}
			The last equality follows from \Cref{lem:projectiveWeight}.
		\end{proof}

		\begin{theorem} \label{thm:xi_densityProjective}
			Let $X \subset \P^{N}$ be an $n$-dimensional projective algebraic manifold defined over $K$. Let $f \colon X \to \R_{\geq 0}$ be a probability density with respect to the volume measure $\mu_X$ on $X$. Let $\bm{\widetilde{A}}$ be the random variable in $K^{n \times N}$ with distribution 
			\[
			    \frac{1 - q^{-(n+1)}}{1 - q^{-1}} \overline{f}(A) 1_{A \in \Ocal^{n \times N}} \mathrm{d}A.
			\]
			Let $\bm{\xi}$ be the random variable obtained by intersecting $X$ with the random space $\Lcal_{\bm{\widetilde{A}}}$ and choosing a point $x$ in the finite set $X \cap \Lcal_{\bm{\widetilde{A}}}$ with probability 
			\[
				 \frac{f(x)}{\overline{f}(\bm{\widetilde{A}})}. 
			\]
			Then $\bm{\xi}$ has density $f$ with respect to $\mu_X$. 
		\end{theorem}
		\begin{proof}
		Let $(\bm{\widetilde{A}}, \bm{\xi})$ be the random variable with values in $\Xcal$ (as defined in \Cref{prop:TangentSpace}) such that $\bm{\widetilde{A}}$ has distribution
		\[
	        \overline{f}(A) \ 1_{A \in \Ocal^{n \times N}} \ \mathrm{d}A
		\] 
		and, given $\bm{\widetilde{A}}$, $\bm{\xi}$ is a random point in $\Lcal_{\bm{\widetilde{A}}} \cap X$ with probability 
		\[
			P(\bm{\xi} = x | \bm{\widetilde{A}}) = \frac{ f(x)}{\overline{f}(\bm{\widetilde{A}})}.
		\]
		Then, by virtue of \Cref{cor:densities}-\ref{cor:densitiespart2} applied to the projection map $\varphi : \Xcal \to K^{n \times N}$, we deduce that the density of $(\bm{\widetilde{A}}, \bm{\xi})$, with respect to $\mu_{\Xcal}$,  is given by
		\begin{align*}
		    f_{\bm{\widetilde{A}}, \bm{\xi}}(A, x) &= \frac{1 - q^{-(n+1)}}{1 - q^{-1}} \overline{f}(A) 1_{A \in \Ocal^{n \times N}} \frac{f(x)}{\overline{f}(A)} \J(\varphi, (A,x)) \\
		    &= \frac{1 - q^{-(n+1)}}{1 - q^{-1}} f(x) \J(\varphi, (A,x)),
		\end{align*}
		for $(A,x) \in \Xcal$. Applying \Cref{prop:TangentSpace} and Corollary \ref{cor:densities}-\ref{cor:densitiespart1} to the map $\phi : \Xcal \to X$, we then deduce that the density of $\bm{\xi}$ is:
		\begin{align*}
			f_{\bm{\xi}}(x) &= \int_{\phi^{-1}(x)} \frac{1 - q^{-(n+1)}}{1 - q^{-1}} f(x) \frac{\J(\varphi, (A,x))}{\J(\phi, (A,x))} \mu_{\phi^{-1}(x)}(\mathrm{d}A)\\
							&= \frac{1 - q^{-(n+1)}}{1 - q^{-1}} f(x) \int_{\phi^{-1}(x)}  \frac{\J(\varphi, (A,x))}{\J(\phi, (A,x))} \mu_{\phi^{-1}(x)}(\mathrm{d}A)\\
							&= \frac{1 - q^{-(n+1)}}{1 - q^{-1}} f(x) \int_{\phi^{-1}(x)}  |\mathrm{det}(A_{|T_x X})| \mu_{\phi^{-1}(x)}(\mathrm{d}A)\\
							&= f(x).
		\end{align*}
		The last equation follows from the \Cref{lem:projectiveWeight}. This concludes the proof.
		\end{proof}

    \section{Sampling linear spaces in practice}\label{sec:Sampling_practically}
		
		In this section we explain how to sample the random planes $\Lcal_{\bm{A}, \bm{b}}$ and $\Lcal_{\bm{A}}$ explicitly. We also explain how to sample the random planes $\Lcal_{\widetilde{\bm{A}}, \widetilde{\bm{b}}}$ from  \Cref{thm:xi_density} and $\Lcal_{\widetilde{\bm{A}}}$ from \Cref{thm:xi_densityProjective} by rejection sampling, and we give bounds on how efficient this sampling method is.
		
		\subsection{Sampling linear spaces explicitly}
		
		When the codimension of the manifold $X$ is small (hypersurfaces for example), for computational reasons, it is easier to find the intersection of $X$ with a linear space $\Ecal$ of complementary dimension $N-n$ when the latter has an explicit form. That is writing $\Ecal$ in the form
		\[
		    \Ecal = u + \mathrm{span}_K(x_1,\dots, x_{N-n}),
		\]
		where $u \in K^N$ and $x_1, \dots, x_{N-n} \in K^N$ are linearly independent.

		\begin{lemma}
		    Let $\bm{A} \in K^{n \times N}$, $\bm{b} \in K^n$ and $\bm{B} \in K^{ (N+1) \times (N-n+1)}$ be matrices with random i.i.d \footnote{Here i.i.d stands for independent and identically distributed.} entries uniformly distributed in $\Ocal$, and $\bm{u}, \bm{x_1}, \dots, \bm{x_{N-n}} \in K^{N}$ be such that
		    $$\begin{pmatrix} \bm{u} \\ 1 \end{pmatrix}, \begin{pmatrix} \bm{x_1} \\ 0 \end{pmatrix}, \dots , \begin{pmatrix} \bm{x_{N-n}} \\ 0 \end{pmatrix} \text{ form an orthonormal basis of } \mathrm{colspan}\left(\bm{B}\right).$$
		    The random affine space $\Ecal_{\bm{u}, \bm{x_1}, \dots, \bm{x_{N-n}}} \coloneqq \bm{u} + \mathrm{span}( \bm{x_1}, \dots, \bm{x_{N-n}})$ has the same probability distribution as $\Lcal_{\bm{A}, \bm{b}}$.
		\end{lemma} 
		\begin{proof}
		    Notice that the linear space $\Lcal_{\bm{A},\bm{b}}$ can be written as
		    \begin{align*}
		        \Lcal_{\bm{A},\bm{b}}  = \left\{ x \in K^{N} \colon(\bm{A} | - \bm{b}) \begin{pmatrix} x \\1 \end{pmatrix} = 0 \right \}.
		    \end{align*}
		    So it suffices to show that $\mathrm{colspan}(\bm{B})$ and $\Ker((\bm{A} | - \bm{b}))$ have the same distribution in the Grassmannian  $\Gr(N - n + 1, K^{N+1})$. Thanks to \Cref{lem:uniqueOrthogonalProb}, it is enough to notice that the distributions of $\mathrm{colspan}(\bm{B})$ and $\Ker((\bm{A} | - \bm{b}))$ are both orthogonally invariant. This is indeed the case since for any $U \in \GL(N+1, \Ocal)$ we have\footnote{ By ``$\stackrel{d}{=}$'' we mean equality in distribution.}:
		    \[
				(\bm{A}|-\bm{b}) U  \stackrel{d}{=} (\bm{A} | -\bm{b}) \quad \text{and} \quad   U  \bm{B} \stackrel{d}{=} \bm{B}.  
			\]
		\end{proof}
	   
		\subsection{Rejection sampling}
			
			Let $X \subset \A^N$ be an affine algebraic manifold of dimension $n$ and degree $d$ and let $f \colon X \to \R_{\geq 0}$ be a probability density function with respect to $\mu_X$.  We recall that the average function $\overline{f}$ in the affine case is defined as:
			\[
				\overline{f}(A,b) = \sum_{ x \in \Lcal_{A,b} \cap X} w_X(x) f(x), \quad \text{ for } (A,b) \in K^{n \times N} \times K^n,
			\]
			where, by convention, the sum is $0$ whenever the intersection $\Lcal_{A,b} \cap X$ is empty or infinite.
			\begin{proposition}[Rejection sampling]
				Suppose that there exists a constant $M > 0$ such that $\overline{f}(A,b) < M$ almost everywhere with respect to $\mathrm{d}A \mathrm{d}b$. 
                Let $(\bm{A}, \bm{b})$ be the random variable with distribution $1_{A \in \Ocal^{n \times N},  b \in \Ocal^n} \mathrm{d}A \mathrm{d}b$ and let $\bm{\eta}$ be a $\{0,1\}$-valued random variable such that:
				$$ P(\bm{\eta} = 1 | (\bm{A},\bm{b}) )	 = \frac{\overline{f}(\bm{A},\bm{b})}{M} \quad  \text{and} \quad	 P(\bm{\eta} = 0 | (\bm{A},\bm{b}) ) = \frac{M - \overline{f}(\bm{A},\bm{b})}{M}.$$
				Then, conditioned on the event $(\bm{\eta} = 1)$, the random variable $(\bm{A}, \bm{b})$ has distribution 
                    $$ \overline{f}(A,b) 1_{A\in \Ocal^{n \times N}, b \in \Ocal^n} \mathrm{d}A \mathrm{d}b. $$
			\end{proposition}
			\begin{proof}
				This follows directly from Bayes' rule as follows
				\begin{align*}
				P\Big((\bm{A},\bm{b}) \in (\mathrm{d}A, \mathrm{d}b) | \eta = 1 \Big) &= \frac{P\Big(\eta = 1 | (\bm{A},\bm{b}) \in (\mathrm{d}A, \mathrm{d}b) \Big) P\Big( (\bm{A},\bm{b}) \in (\mathrm{d}A,\mathrm{d}b) \Big) }{P(\bm{\eta} = 1)} \\
				&= \frac{P\Big(\eta = 1 | (\bm{A},\bm{b}) \in (\mathrm{d}A,\mathrm{d}b)\Big) }{P(\bm{\eta} = 1)} 1_{A \in \Ocal^{n \times N}, \  b \in \Ocal^n} \ \mathrm{d}A \mathrm{d}b\\
				&= \frac{  \overline{f}(A,b) / M }{\E[\overline{f}(\bm{A}, \bm{b})] / M} 1_{A \in \Ocal^{n \times N}, \ b \in \Ocal^n} \ \mathrm{d}A \mathrm{d}b \\
				& = \overline{f}(A,b).
				\end{align*}
				The last equation follows from \Cref{thm:integral_EP} and the equality 
				\[
    				P(\bm{\eta} = 1) = \E[P(\bm{\eta} = 1 | (\bm{A}, \bm{b}))] = \frac{1}{M} \E[\overline{f}(\bm{A},\bm{b})] = \frac{1}{M}.  
				\]
				
			\end{proof}
		
			\begin{lemma}\label{lem:AffineBounds}
				Let $f \colon X \to \R_{\geq 0}$ be a probability density function supported on $X \cap \varpi^{-r} \Ocal^N$ for some integer $r \geq 0$. Suppose that $\kappa \coloneqq \sup_{x \in X } f(x) < \infty$. Then we have
				\[
					\overline{f}(A,b) \leq d q^{(n+1)r} \ \frac{1 - q^{-(n+1)}}{1 - q^{-1}} \kappa,
				\]
                with $d \coloneqq \deg(X)$. In particular, if $f$ is the uniform probability density on $X \cap \Ocal^N$ then:
				\[
					\overline{f}(A,b) \leq d \ \frac{1 - q^{-(n+1)}}{1 - q^{-1}}.
				\]
			\end{lemma}
			\begin{proof}
			    Let $x \in X$ and $U,W,S_x$ as in \Cref{def:NrDefAffine}. Then, since the columns of $UW$ are orthonormal in $K^N$, its rows are in $\Ocal^n$ and, modulo $\varpi$, they span $k^n$. So we deduce that
			    \[
			    \Nr(X,x) = \N(S_x U W) \geq \min(1,\norm{x}^{-1}).
			    \]
			    Hence, from \Cref{def:affineWeight}, we get 
			    \[
			    w_X(X,x) \leq \frac{1 - q^{-(n+1)}}{1 - q^{-1}} \max(1, \norm{x}^{n+1}).
			    \]
			    Then for $(A,b) \in K^{n \times N} \times K^n$ we get
			    \[
			       \overline{f}(A,b) \leq \#(X \cap \Lcal_{A,b}) \frac{1 - q^{-(n+1)}}{1 - q^{-1}} q^{(n+1)r} \sup_{x \in X}f(x).
			    \]
			    Since the number of intersection points $\#(X \cap \Lcal_{A,b})$ is at most $d = \deg(X)$ (except for a measure zero set of $(A,b) \in K^{n \times N} \times K^n$, see \Cref{subsec:intersection}), we deduce the desired result. The second statement is an immediate consequence of the first one.
			\end{proof}
			
			\begin{remark}
			    The bound given for $\overline{f}(A,b)$ is far from being sharp. Moreover, when one wishes to sample from $X \cap \varpi^{-r} \Ocal^N$, this bound is not very practical for rejection sampling. In this case, it is better to use \Cref{rem:weightBehaviour}~\ref{rem:tip}.
			\end{remark}

		    Let $h \colon X \to \R$ be an integrable function on $X$ supported on $X \cap \Ocal^{N}$ and let $(\bm{A}_i, \bm{b}_i)_{i \geq 0}$ be a sequence of i.i.d random variables such that $(\bm{A}_i,\bm{b}_i)$ has the uniform distribution on $\Ocal^{n  \times N} \times \Ocal^n$ for all $i \geq 0$. Finally, set: 
		    \[
		    S_m(h) \coloneqq \overline{h}(\bm{A}_1, \bm{b}_1) + \overline{h}(\bm{A}_2, \bm{b}_2) + \cdots + \overline{h}(\bm{A}_m, \bm{b}_m).
		    \]
		    Then we have the following:
		    \begin{proposition}
		        The random variable $S_m(h) / m$ converges almost surely to the integral $I(h) \coloneqq \int_X h(x)\mu_X(\mathrm{d}x)$ as $m \uparrow \infty$. Moreover, if $ \kappa \coloneqq \sup_{x \in X} |h(x)| < \infty$, then
		        \[
		        P \left( \left| \frac{S_m(h)}{m}  - I(h) \right| \geq \epsilon \right) \leq  \frac{ \kappa^2 d^2}{\epsilon^2 m} \left(\frac{1 - q^{-(n+1)}}{1 - q^{-1}} \right)^2 \quad \text{for } m \geq 1,
		        \]
                with $d \coloneqq \deg(X)$.
		    \end{proposition}
		    \begin{proof}
		        The first statement is an immediate application of the law of large numbers. The second follows from \Cref{lem:AffineBounds} and Chebychev's inequality.
		    \end{proof}
		    
		    \begin{remark}
		        While this section focuses on affine manifolds, the results discussed within can be restated and proved for projective manifolds without much difficulty.
		    \end{remark}
		  
\section{Applications and examples} \label{sec:Sampling_applications}
        
            In this section we discuss a few concrete examples and applications. The first case of interest is when the algebraic manifold $X$ is an algebraic group.
            
    		\subsection{Measures on algebraic groups}
    		
    		Let $G$ be an algebraic group defined over $K$, by which we mean a smooth (either affine or projective) algebraic variety together with
    		
    		\begin{enumerate}
    		\item (\emph{identity element}) an element $e \in G$,
            \item (\emph{multiplication}) a  morphism $m \colon G \times G \to G, (x,y) \mapsto xy$,
            \item (\emph{inverse}) a morphism $\iota \colon G \to G, x \mapsto x^{-1}$,
    		\end{enumerate} 
    		with respect to which $G$ is a group (see \cite{MilneAlgGroups} or \cite{BorelAlgGroups} for a detailed account). In our discussion, $m$ and $\iota$ are $K$-morphisms and we are interested in the group $G(K)$ of $K$ points of $G$ which we also denote by $G$ for simplicity and, for our purposes, $G$ is embedded in some affine or projective space over $K$.

    		The group $G$ is a locally compact topological group and thus admits a left Haar measure; that is a non-zero measure $\nu_G$ such that
    		\[
    		    \nu(g A) = \nu(A), \quad \text{for any Borel measurable set } A \subset G.
    		\]
    		which is unique up to scaling. If $G$ is an algebraic group embedded in a projective space as an algebraic manifold, then $G$ is compact and the measure $\mu_G$ is then finite and also right-invariant. In this case we normalize $\nu$ so that $\nu(G) = 1$. In the case where $G$ is affine, the measure $\nu$ is finite on the set $G(\Ocal)$ of $\Ocal$-points of $G$ and we normalise $\nu$ so that $\nu(G(\Ocal)) = 1$.

            \begin{remark}
                It is not always the case that the points in $G(\Ocal)$ form a subgroup of $G$. For example, this fails to be the case for $G = \GL(n,K)$.
            \end{remark}
    		
          \begin{example}
        		Let $n \geq 1$ be a positive integer. If $G$ is either the special linear group $\SL(n,K)$ or the special orthogonal group $\SO(n,K)$ or the symplectic group $\Sp(n,K)$, the $\Ocal$-points $G(\Ocal)$ form a compact subgroup of $G$. Moreover, the normalized Haar measure $\nu$ on $G(\Ocal)$ coincides with the uniform probability  measure on $G(\Ocal)$ with respect to the volume measure $\mu_G$ as defined in \Cref{eq:defOfAffineMeasure} (or in \Cref{eq:defOfProjMeasure} for the projective case). This is because the measure $\mu_{\A^{n \times n}}$ is invariant under the action of $\GL(n , \Ocal)$ and in particular under the action of $G(\Ocal) \subset \GL(n ,\Ocal)$ and hence $\mu_G$ is also $G(\Ocal)$-invariant. So using \Cref{thm:xi_density} we can sample from the Haar measure on the compact matrix groups $\SL(n,\Ocal)$, $\SO(n,\Ocal)$ and $\Sp(n,\Ocal)$. For small values of $n$, we provide examples of this in the repository \cite{MathRepo}.
    		\end{example}

    		\subsection{Moduli spaces}
    		    
    		    Another case of interest is when the algebraic manifold $X$ is a moduli space parameterizing certain objects. Then sampling from $X$, we can get an idea of how often a certain property of these objects holds or how rare are objects of certain kind are in $X$. Here, we give two examples of such a situation.
    		    
    		\subsubsection{Modular curves}
    		
        		Let $N$ be a positive integer and consider the modular curve $X_1(N)$. This is a smooth projective curve defined over $\Q$, and it has the following \textit{moduli} interpretation: for any field $K$ with characteristic $0$, non-cuspidal\footnote{Modular curves have only finitely many cuspidal points. This will be important for what follows.} $K$-points of $X_1(N)$ parametrize isomorphism classes of pairs $(E,P)$, where $E$ is an elliptic curve over $K$ and $P$ is a point of $E(K)$ of order $N$. For a reference, see \cite{DiamondShurman}, or \cite[Section~C.13]{SilvermanI} for a quick introduction.

                \medskip 
                
        		In this example, we will sample uniformly from $\Z_{31}$-points of $X_1(30)$, and compute the \textit{Tamagawa numbers} of the corresponding elliptic curves over $\Q_{31}$. For an elliptic curve $E/\Qp$, the finite index
        		\[c_p = [E(\Qp) : E^0(\Qp)]\]
        		is referred to as the Tamagawa number of $E/\Qp$, where $E^0(\Qp)$ is the subgroup of $E(\Qp)$ consisting of points that have good reduction. Clearly, if $E/\Qp$ has good reduction, then $c_p$ equals $1$. We note that Tamagawa numbers of elliptic curves are important local arithmetic invariants. They arise in the conjecture of Birch and Swinnerton-Dyer, see \cite[Section~C.16]{SilvermanI}. Moreover, they can be easily computed using Magma \cite{magma}.

            \medskip
        
        		The following (optimized) equation for $X_1(30)$ was provided by Sutherland in \cite{Sutherland12}:
        		\begin{align*} 
                X_1(30) : y^6 &+ (x^6 - 5x^5 + 6x^4 + 3x^3 - 6x^2 + 7x + 3)y^5 \\
                        &+ (x^7 - 3x^6 - 13x^5 + 44x^4 - 18x^3 + x^2 + 18x + 3)y^4 \\
                        &+ (x^8 - 3x^7 - 13x^6 + 27x^5 + 46x^4 - 32x^3 + 21x^2 + 15x + 1)y^3 \\
                        &+ 2x(x^7 - 8x^6 + 9x^5 + 20x^4 + 6x^3 - 6x^2 + 9x + 2)y^2 \\
                        &- 4x^2(2x^5 - 7x^4 - 3x^3 - 1)y + 8x^6 = 0.
                \end{align*}
                Moreover, if $(x_0,y_0)$ is a non-cuspidal point on $X_1(30)$, then the corresponding elliptic curve is of the form
                \[
                    y^2 = x^3 + (t^2-2qt-2)x^2 -(t^2-1)(qt+1)^2x,
                \]
                where
                \begin{align*}
                    q &= y_0+1, \\
                    t &= 4(y_0+1)(x_0+y_0)/(x_0y_0^3 - 4x_0y_0 - 4x_0 - 3y_0^3 - 6y_0^2 - 4y_0).
                \end{align*}
                See the table in \text{\url{https://math.mit.edu/~drew/X1_optcurves.html}}. Table~\ref{tab:TamNum} presents the Tamagawa numbers of elliptic curves obtained for a sample of $500000$ $\Z_{31}$-points on $X_1(30)$, and the number of times they occurred.
        
        		\begin{table}[H]
        		    \centering
        	    	\begin{tabular}{|c|c||c|c||c|c|}	   
        				\hline
        				 $c_{31}$ &      & $c_{31}$ &  & $c_{31}$ &   \\[0.6ex]
        	    		\hline
        			        1    & 267228  & 8 & 2 & 20 & 380 \\[0.6ex]
        				\hline	
        				    2    &  52832 & 9 & 69 & 24 & 3 \\[0.6ex]
        				\hline
        			        3    & 56509 & 10 & 13163 & 30 & 6520 \\[0.6ex]
        				\hline
        					4	& 1669 & 12 & 1781 & 45 & 10 \\[0.6ex]
        	    		\hline
        			        5    & 27673  & 15 & 12980 & 60 & 180 \\[0.6ex]
        				\hline	
        				    6   &  58927  & 18 & 69 & 90 & 5 \\[0.6ex]
        				\hline
        			\end{tabular}
        		\caption{The Tamagawa numbers and their multiplicities that appeared in our sampling.}
        		\label{tab:TamNum}		
        		\end{table}
    		
        \subsubsection{Hilbert modular surfaces}

        Following the notation in \cite{ElkiesKumar}, we will work with \textit{Hilbert modular surfaces} $Y_{-}(D)$. These surfaces parametrize abelian surfaces with real multiplication. More precisely, let $d>1$ be a square-free integer, and set
        \[ D =
            \begin{cases}
                \hfil d       & \text{if } d\equiv 1 \text{ mod } 4,\\
                4d  & \text{if } d\equiv 2,3 \text{ mod } 4.
            \end{cases} \]
        Note that $D$ is nothing but the discriminant of the ring of integers $\Ocal_D$ of the real quadratic field $\Q(\sqrt{D})$. Such a number is called a \textit{positive fundamental discriminant}. The quotient 
        \[\text{PSL}_2(\Ocal_D) \setminus \left( \Hcal^{+}\times\Hcal^{-} \right)\]
        is the coarse moduli space of principally polarized abelian surfaces with real multiplication by $\Ocal_D$. Here, $\Hcal^{+}$ (resp. $\Hcal^{-}$) denotes the complex upper (resp. lower) half plane. There is a holomorphic map from this quotient to the moduli space $\Acal_2$ of principally polarized abelian surfaces. The image is the \textit{Humbert surface} $\Hcal_D$, and the Hilbert modular surface $Y_{-}(D)$ is a double cover of $\Hcal_D$ branched along a finite union of modular curves. For the theory of Hilbert modular surfaces, see, for example,  \cite{vanderGeer,Bruinier}.
    
        The surfaces $Y_{-}(D)$ have models over $\Q$, and points on these surfaces correspond, generically, to Jacobians of smooth projective curves\footnote{Recall that a principally polarized abelian surface over an algebraically closed field is either the Jacobian variety of a smooth projective curve of genus $2$ or the product of two elliptic curves.} of genus $2$. Explicit equations for birational models of $Y_{-}(D)$, as well as the Igusa--Clebsch invariants $I_2,I_4,I_6$ and  $I_{10}$ of the corresponding genus-$2$ curves, were provided by Elkies and Kumar in \cite{ElkiesKumar} for all fundamental discriminants $D$ between $1$ and $100$. In this final example, we will 
        \begin{itemize}
            \item sample uniformly from $\Z_5$-points of $Y_{-}(5)$, and
            \item compute the minimal skeleta of the Berkovich analytification of the corresponding genus-$2$ curves. 
        \end{itemize}
        It is well know that there are precisely $7$ different (graph-theoretical) types, which are depicted in Figure~\ref{GraphsOfMinSke}. The recent work of Helminck \cite{helminck2021Igusa} shows that \textit{tropical Igusa invariants}, which can easily be computed from Igusa--Clebsch invariants, distinguish between the different types; see \cite[Theorem~2.11]{helminck2021Igusa}. See also \cite{HelminckElMaazouzKaya} for a similar result concerning Picard curves.
        
        \begin{figure}[ht]
    	    \centering
         \scalebox{0.8}{
        	    \begin{tikzpicture}
        	    
                    \filldraw (0.75,0) circle (3.5pt);
                    \filldraw (0.75,0.4) node{$2$};                    
                    \filldraw (0.75,-1) node{\textbf{Type~$\mathrm{I}$}};
                    \qquad
                    \filldraw (4,0) circle (3.5pt);
                    \filldraw (4.35,0) node{$1$};
                    \draw[very thick] (3.25,0) ellipse (0.75 and 0.42);
                    \filldraw (3.35,-1) node{\textbf{Type~$\mathrm{II}$}};

                    \qquad
                    
                    \filldraw (7,0) circle (3.5pt);
                    \draw[very thick] (6.25,0) ellipse (0.75 and 0.42);
                    \draw[very thick] (7.75,0) ellipse (0.75 and 0.42);
                    
                    \filldraw (7,-1) node{\textbf{Type~$\mathrm{III}$}};
                    
                    \qquad
                    
                    \filldraw (10,0) circle (3.5pt);
                    \filldraw (11.5,0) circle (3.5pt);
                    \draw [very thick] (10,0) -- (11.5,0);
                    \draw[very thick] (10.75,0) ellipse (0.75 and 0.42);
                    
                    \filldraw (10.75,-1) node{\textbf{Type~$\mathrm{IV}$}};

                    \filldraw (0,-2.5) circle (3.5pt);
                    \filldraw (1.5,-2.5) circle (3.5pt);
                    \filldraw (0,-2.1) node{$1$};
                    \filldraw (1.5,-2.1) node{$1$};
                    \draw [very thick] (0,-2.5) -- (1.5,-2.5);
                    \filldraw (0.75,-3.5) node{\textbf{Type~$\mathrm{V}$}};
                    
                    \qquad

                    \filldraw (5,-2.5) circle (3.5pt);
                    \filldraw (6.5,-2.5) circle (3.5pt);
                    \filldraw (6.85,-2.5) node{$1$};
                    \draw[very thick] (4.25,-2.5) ellipse (0.75 and 0.42);
                    \draw [very thick] (5,-2.5) -- (6.5,-2.5);
                    
                    \filldraw (5.2,-3.5) node{\textbf{Type~$\mathrm{VI}$}};
                    
                    \qquad
                    
                    \filldraw (10,-2.5) circle (3.5pt);
                    \filldraw (11.5,-2.5) circle (3.5pt);
                    \draw[very thick] (9.25,-2.5) ellipse (0.75 and 0.42);
                    \draw [very thick] (10,-2.5) -- (11.5,-2.5);
                    \draw[very thick] (12.25,-2.5) ellipse (0.75 and 0.42);
                    \filldraw (10.75,-3.5) node{\textbf{Type~$\mathrm{VII}$}};
                    
                    \end{tikzpicture}}
    	    \caption{Minimal skeleta of the Berkovich analytification of genus-$2$ curves.}
    	    \label{GraphsOfMinSke}
    	\end{figure}
        
        A birational model of the surface $Y_{-}(5)$ is given by
        \[
            z^2 = 2\left( 6250h^2 - 4500g^2h - 1350gh - 108h - 972g^5 -324g^4 - 27g^3 \right),
        \]
        see \cite[Theorem~16]{ElkiesKumar}. Moreover, the map from from $\Hcal_{5}$ to $\Acal_2$ (or, more precisely, to the moduli space $\Mcal_2$ of curves of genus $2$) is given by
        \[\left( I_2 :I_4 : I_6 : I_{10} \right) = \left( 6(4g+1), 9g^2, 9(4h + 9g^3 + 2g^2), 4h^2 \right)\]
        see \cite[Corollary~15]{ElkiesKumar}. Table~\ref{tab:Types} shows how many times the types occurred for a sample of $500000$ $\Z_{5}$-points on $Y_{-}(5)$. 
        
        \begin{table}[H]
		    \centering
	    	\begin{tabular}{|c|c|c|c|c|c|c|}	   
				\hline
				 \textbf{Type~$\mathrm{I}$} & \textbf{Type~$\mathrm{II}$} &  \textbf{Type~$\mathrm{III}$} & \textbf{Type~$\mathrm{IV}$} & \textbf{Type~$\mathrm{V}$} & \textbf{Type~$\mathrm{VI}$} & \textbf{Type~$\mathrm{VII}$}  \\[0.6ex]
	    		\hline
			       415075 & 0 & 40108 & 23166 & 21624 & 1 & 26 \\[0.6ex]
				\hline
			\end{tabular}
		    \caption{The multiplicities of the types that appeared in our sampling.}
		    \label{tab:Types}		
    	\end{table} 
    	
    	As shown the in table, Types $\mathrm{II}$, $\mathrm{VI}$ and $\mathrm{VII}$ are quite rare. In fact, we never see Type~$\mathrm{II}$, and it is unclear to the authors if there is a theoretical reason behind this.

\section*{Acknowledgments}
This work was initiated during YE's visit to the Max Planck Institute for Mathematics in the Sciences. We are grateful to the institute for its warm hospitality and inspiring research environment. We extend our thanks to Peter B\"urgisser, Avinash Kulkarni, and Antonio Lerario for kindly sharing an early draft of their manuscript \cite{BKL22}. Special thanks go to Corrine Elliott, Steven N. Evans, Kiran Kedlaya, Bartosz Naskrecki, Bernd Sturmfels, and Tristan Vaccon for their insightful discussions and valuable feedback. Finally, we thank the anonymous referees for their thoughtful suggestions, which greatly improved the presentation and readability of this manuscript.


\begin{thebibliography}{10}

    \bibitem{manssour2020probabilistic}
    {\sc Ait El~Manssour, R., and Lerario, A.}
    \newblock Probabilistic enumerative geometry over {$p$}-adic numbers: linear
      spaces on complete intersections.
    \newblock {\em Ann. H. Lebesgue 5\/} (2022), 1329--1360.
    
    \bibitem{andersenDiac}
    {\sc Andersen, H.~C., and Diaconis, P.}
    \newblock Hit and run as a unifying device.
    \newblock {\em J. Soc. Fr. Stat. \& Rev. Stat. Appl. 148}, 4 (2007), 5--28.
    
    \bibitem{BhargavaCremonaFisherGajovic21}
    {\sc Bhargava, M., Cremona, J., Fisher, T., and Gajovi{\'c}, S.}
    \newblock The density of polynomials of degree $n$ over $\mathbb{Z}_p$ having
      exactly $r$ roots in $\mathbb{Q}_p$.
    \newblock {\em Proceedings of the London Mathematical Society\/} (2022).
    
    \bibitem{BorelAlgGroups}
    {\sc Borel, A.}
    \newblock {\em Linear algebraic groups}, second~ed., vol.~126 of {\em Graduate
      Texts in Mathematics}.
    \newblock Springer-Verlag, New York, 1991.
    
    \bibitem{magma}
    {\sc Bosma, W., Cannon, J., and Playoust, C.}
    \newblock The {M}agma algebra system. {I}. {T}he user language.
    \newblock {\em J. Symbolic Comput. 24}, 3-4 (1997), 235--265.
    \newblock Computational algebra and number theory (London, 1993).
    
    \bibitem{Bostan2005}
    {\sc Bostan, A., González-Vega, L., Perdry, H., and Éric Schost}.
    \newblock From newton sums to coefficients: complexity issues in characteristic
      \( p \).
    \newblock In {\em Proceedings of MEGA'05\/} (2005).
    
    \bibitem{RandomPoints}
    {\sc Breiding, P., and Marigliano, O.}
    \newblock Random points on an algebraic manifold.
    \newblock {\em SIAM J. Math. Data Sci. 2}, 3 (2020), 683--704.
    
    \bibitem{Bruinier}
    {\sc Bruinier, J.~H.}
    \newblock Hilbert modular forms and their applications.
    \newblock In {\em The 1-2-3 of modular forms}, Universitext. Springer, Berlin,
      2008, pp.~105--179.
    
    \bibitem{BKL22}
    {\sc B{\"u}rgisser, P., Kulkarni, A., and Lerario, A.}
    \newblock Nonarchimedean integral geometry.
    \newblock {\em arXiv:2206.03708\/} (2023).
    
    \bibitem{Caruso2017}
    {\sc Caruso, X.}
    \newblock Computations with \( p \)-adic numbers.
    \newblock {\em Les cours du CIRM 5}, 1 (2017), 75.
    
    \bibitem{caruso2021zeroes}
    {\sc Caruso, X.}
    \newblock Where are the zeroes of a random {$p$}-adic polynomial?
    \newblock {\em Forum Math. Sigma 10\/} (2022), Paper No. e55, 41.
    
    \bibitem{DiacanAl}
    {\sc Diaconis, P., Lebeau, G., and Michel, L.}
    \newblock Geometric analysis for the {M}etropolis algorithm on {L}ipschitz
      domains.
    \newblock {\em Invent. Math. 185}, 2 (2011), 239--281.
    
    \bibitem{saloffCosteDiac}
    {\sc Diaconis, P., and Saloff-Coste, L.}
    \newblock What do we know about the {M}etropolis algorithm?
    \newblock {\em J. Comput. System Sci. 57}, 1 (1998), 20--36.
    \newblock 27th Annual ACM Symposium on the Theory of Computing (STOC'95) (Las
      Vegas, NV).
    
    \bibitem{DiamondShurman}
    {\sc Diamond, F., and Shurman, J.}
    \newblock {\em A first course in modular forms}, vol.~228 of {\em Graduate
      Texts in Mathematics}.
    \newblock Springer-Verlag, New York, 2005.
    
    \bibitem{3264Eisenbud}
    {\sc Eisenbud, D., and Harris, J.}
    \newblock {\em 3264 and all that---a second course in algebraic geometry}.
    \newblock Cambridge University Press, Cambridge, 2016.
    
    \bibitem{maazouz2021gaussian}
    {\sc El~Maazouz, Y.}
    \newblock The {G}aussian entropy map in valued fields.
    \newblock {\em Algebr. Stat. 13}, 1 (2022), 1--18.
    
    \bibitem{MathRepo}
    {\sc {El Maazouz}, Y., and {Kaya}, E.}
    \newblock Mathrepo code repository.
    \newblock Available at {\tt
      \url{https://mathrepo.mis.mpg.de/SamplingpAdicManifolds/index.html}}.
    
    \bibitem{maazouz2019statistics}
    {\sc El~Maazouz, Y., and Tran, N.~M.}
    \newblock Statistics and tropicalization of local field {G}aussian measures.
    \newblock {\em arXiv:1909.00559\/} (2019).
    
    \bibitem{ElkiesKumar}
    {\sc Elkies, N., and Kumar, A.}
    \newblock K3 surfaces and equations for {H}ilbert modular surfaces.
    \newblock {\em Algebra Number Theory 8}, 10 (2014), 2297--2411.
    
    \bibitem{evansBMLocalField}
    {\sc Evans, S.~N.}
    \newblock Local field {B}rownian motion.
    \newblock {\em J. Theoret. Probab. 6}, 4 (1993), 817--850.
    
    \bibitem{evansNoise95}
    {\sc Evans, S.~N.}
    \newblock {$p$}-adic white noise, chaos expansions, and stochastic integration.
    \newblock In {\em Probability measures on groups and related structures, {XI}
      ({O}berwolfach, 1994)}. World Sci. Publ., River Edge, NJ, 1995, pp.~102--115.
    
    \bibitem{evans2001local}
    {\sc Evans, S.~N.}
    \newblock Local fields, {G}aussian measures, and {B}rownian motions.
    \newblock In {\em Topics in probability and {L}ie groups: boundary theory},
      vol.~28 of {\em CRM Proc. Lecture Notes}. Amer. Math. Soc., Providence, RI,
      2001, pp.~11--50.
    
    \bibitem{ElementaryDivisors}
    {\sc Evans, S.~N.}
    \newblock Elementary divisors and determinants of random matrices over a local
      field.
    \newblock {\em Stochastic Process. Appl. 102}, 1 (2002), 89--102.
    
    \bibitem{GHJMPSST}
    {\sc Gubser, S.~S., Heydeman, M., Jepsen, C., Marcolli, M., Parikh, S., Saberi,
      I., Stoica, B., and Trundy, B.}
    \newblock Edge length dynamics on graphs with applications to {$p$}-adic
      {A}d{S}/{CFT}.
    \newblock {\em J. High Energy Phys.}, 6 (2017), 157, front matter+34.
    
    \bibitem{pAdicAdSCFT}
    {\sc Gubser, S.~S., Knaute, J., Parikh, S., Samberg, A., and Witaszczyk, P.}
    \newblock {$p$}-adic {A}d{S}/{CFT}.
    \newblock {\em Comm. Math. Phys. 352}, 3 (2017), 1019--1059.
    
    \bibitem{harrisFirstCourse}
    {\sc Harris, J.}
    \newblock {\em Algebraic geometry}, vol.~133 of {\em Graduate Texts in
      Mathematics}.
    \newblock Springer-Verlag, New York, 1995.
    \newblock A first course, Corrected reprint of the 1992 original.
    
    \bibitem{Hasse1}
    {\sc Hasse, H.}
    \newblock \"{U}ber die {D}arstellbarkeit von {Z}ahlen durch quadratische
      {F}ormen im {K}\"{o}rper der rationalen {Z}ahlen.
    \newblock {\em J. Reine Angew. Math. 152\/} (1923), 129--148.
    
    \bibitem{Hasse2}
    {\sc Hasse, H.}
    \newblock Darstellbarkeit von zahlen durch quadratische formen in einem
      beliebigen algebraischen zahlkörper.
    \newblock {\em Journal für die reine und angewandte Mathematik 153\/} (1924),
      113--130.
    
    \bibitem{helminck2021Igusa}
    {\sc Helminck, P.~A.}
    \newblock Tropical {I}gusa {I}nvariants.
    \newblock {\em arXiv 1604.03987\/} (2021).
    
    \bibitem{HelminckElMaazouzKaya}
    {\sc Helminck, P.~A., El~Maazouz, Y., and Kaya, E.}
    \newblock Tropical invariants for binary quintics and reduction types of
      {P}icard curves.
    \newblock {\em Glasg. Math. J. 66}, 1 (2024), 65--87.
    
    \bibitem{Hensel05}
    {\sc Hensel, K.}
    \newblock \"{U}ber eine neue {B}egr\"{u}ndung der {T}heorie der algebraischen
      {Z}ahlen.
    \newblock {\em Jahresbericht der Deutschen Mathematiker- Vereinigung 6\/}
      (1897), 83--88.
    
    \bibitem{HMPS}
    {\sc Heydeman, M., Marcolli, M., Parikh, S., and Saberi, I.}
    \newblock Nonarchimedean holographic entropy from networks of perfect tensors.
    \newblock {\em Adv. Theor. Math. Phys. 25}, 3 (2021), 591--721.
    
    \bibitem{Kedlaya}
    {\sc Kedlaya, K.~S.}
    \newblock Counting points on hyperelliptic curves using {M}onsky-{W}ashnitzer
      cohomology.
    \newblock {\em J. Ramanujan Math. Soc. 16}, 4 (2001), 323--338.
    
    \bibitem{kulkarni2021integGeo}
    {\sc Kulkarni, A., and Lerario, A.}
    \newblock {$p$}-adic integral geometry.
    \newblock {\em SIAM J. Appl. Algebra Geom. 5}, 1 (2021), 28--59.
    
    \bibitem{LLL}
    {\sc Lenstra, A.~K., Lenstra, Jr., H.~W., and Lov\'{a}sz, L.}
    \newblock Factoring polynomials with rational coefficients.
    \newblock {\em Math. Ann. 261}, 4 (1982), 515--534.
    
    \bibitem{Liu}
    {\sc Liu, J.~S.}
    \newblock {\em Monte {C}arlo strategies in scientific computing}.
    \newblock Springer Series in Statistics. Springer-Verlag, New York, 2001.
    
    \bibitem{MilneAlgGroups}
    {\sc Milne, J.~S.}
    \newblock {\em Algebraic groups}, vol.~170 of {\em Cambridge Studies in
      Advanced Mathematics}.
    \newblock Cambridge University Press, Cambridge, 2017.
    \newblock The theory of group schemes of finite type over a field.
    
    \bibitem{Mikowski}
    {\sc Minkowski, H.}
    \newblock {\em Geometrie der {Z}ahlen}.
    \newblock Bibliotheca Mathematica Teubneriana, Band 40. Johnson Reprint Corp.,
      New York-London, 1968.
    
    \bibitem{NeukirchBook}
    {\sc Neukirch, J.}
    \newblock {\em Algebraic number theory}, vol.~322 of {\em Grundlehren der
      mathematischen Wissenschaften [Fundamental Principles of Mathematical
      Sciences]}.
    \newblock Springer-Verlag, Berlin, 1999.
    \newblock Translated from the 1992 German original and with a note by Norbert
      Schappacher, With a foreword by G. Harder.
    
    \bibitem{mpopaInteg}
    {\sc Popa, M.}
    \newblock Chapter~3: $p$-adic integration.
    \newblock \url{https://people.math.harvard.edu/~mpopa/571/}.
    
    \bibitem{LieGroups}
    {\sc Schneider, P.}
    \newblock {\em {$p$}-adic {L}ie groups}, vol.~344 of {\em Grundlehren der
      mathematischen Wissenschaften [Fundamental Principles of Mathematical
      Sciences]}.
    \newblock Springer, Heidelberg, 2011.
    
    \bibitem{SerreLocalField}
    {\sc Serre, J.-P.}
    \newblock {\em Local fields}, vol.~67 of {\em Graduate Texts in Mathematics}.
    \newblock Springer-Verlag, New York-Berlin, 1979.
    \newblock Translated from the French by Marvin Jay Greenberg.
    
    \bibitem{Serre81}
    {\sc Serre, J.-P.}
    \newblock Quelques applications du th\'{e}or\`eme de densit\'{e} de
      {C}hebotarev.
    \newblock {\em Inst. Hautes \'{E}tudes Sci. Publ. Math. 5\/} (1981), 323--401.
    
    \bibitem{SilvermanI}
    {\sc Silverman, J.~H.}
    \newblock {\em The arithmetic of elliptic curves}, second~ed., vol.~106 of {\em
      Graduate Texts in Mathematics}.
    \newblock Springer, Dordrecht, 2009.
    
    \bibitem{Sutherland12}
    {\sc Sutherland, A.~V.}
    \newblock Constructing elliptic curves over finite fields with prescribed
      torsion.
    \newblock {\em Math. Comp. 81}, 278 (2012), 1131--1147.
    
    \bibitem{sagemath}
    {\sc {The Sage Developers}}.
    \newblock {S}age{M}ath, the {S}age {M}athematics {S}oftware {S}ystem ({V}ersion
      9.5), 2022.
    \newblock {\tt https://www.sagemath.org}.
    
    \bibitem{vanderGeer}
    {\sc van~der Geer, G.}
    \newblock {\em Hilbert modular surfaces}, vol.~16 of {\em Ergebnisse der
      Mathematik und ihrer Grenzgebiete (3) [Results in Mathematics and Related
      Areas (3)]}.
    \newblock Springer-Verlag, Berlin, 1988.
    
    \bibitem{vanRooij}
    {\sc van Rooij, A. C.~M.}
    \newblock {\em Non-{A}rchimedean functional analysis}, vol.~51 of {\em
      Monographs and Textbooks in Pure and Applied Mathematics}.
    \newblock Marcel Dekker, Inc., New York, 1978.
    
    \bibitem{padicPhys1}
    {\sc Vladimirov, V.~S., Volovich, I.~V., and Zelenov, E.~I.}
    \newblock {\em {$p$}-adic analysis and mathematical physics}, vol.~1 of {\em
      Series on Soviet and East European Mathematics}.
    \newblock World Scientific Publishing Co., Inc., River Edge, NJ, 1994.
    
    \end{thebibliography}
\end{document}